\documentclass[12pt,a4paper]{article}
\usepackage{latexsym,amssymb,amsmath,mathrsfs}

\usepackage{sectsty}
\usepackage{color}
\definecolor{darkyellow2}{rgb}{0.6,0.6,0.2}
\definecolor{darkgreen}{rgb}{0,0.692157,0}
\definecolor{navyblue}{cmyk}{0.5,0.5,0,0.5}
\usepackage[anchorcolor=black,%
linkcolor=black,%
colorlinks=true,%
citecolor=navyblue,%
bookmarks=true,%
bookmarksnumbered=true,%
]{hyperref}
\usepackage{graphicx} 
\newtheorem{Thm}{Theorem}[section]
\newtheorem{Prop}[Thm]{Proposition}
\newtheorem{Cor}[Thm]{Corollary}
\newtheorem{Lem}[Thm]{Lemma}

\newtheorem{Rem}[Thm]{Remark}
\newenvironment{Proof}{\par\begin{trivlist}%
\item[]{\bf Proof.\ }}%
{\hfill $\square$ \end{trivlist}\par}
\newenvironment{tProof}[1]{\par\begin{trivlist}%
\item[]{\bf Proof of #1.\ }}%
{\hfill $\square$ \end{trivlist}\par}

\makeatletter \@addtoreset{equation}{section} \makeatother

\renewcommand{\P}{\mathbb{P}}
\newcommand{\E}{\mathbb{E}}
\newcommand{\R}{\mathbb{R}}
\newcommand{\N}{\mathbb{N}}

\renewcommand{\a}{\alpha}
\renewcommand{\b}{\beta}
\newcommand{\gm}{\gamma}
\newcommand{\dl}{\delta}
\newcommand{\ep}{\varepsilon}
\newcommand{\te}{\theta}
\newcommand{\lm}{\lambda}
\newcommand{\sg}{\sigma}
\newcommand{\ph}{\varphi}
\newcommand{\Gm}{\Gamma}
\newcommand{\Dl}{\Delta}
\newcommand{\Lm}{\Lambda}
\newcommand{\sL}{\mathscr{L}}
\newcommand{\sP}{\mathscr{P}}

\newcommand{\abs}[1]{\left| #1 \right|}
\newcommand{\abra}[1]{\left( #1 \right)}
\newcommand{\bbra}[1]{\left\{ #1 \right\}}
\newcommand{\cbra}[1]{\left[ #1 \right]}
\newcommand{\dbra}[1]{\langle #1 \rangle}
\DeclareMathOperator{\Ric}{Ric}
\DeclareMathOperator{\diam}{diam}

\DeclareMathOperator{\Cut}{Cut}
\DeclareMathOperator{\Hess}{Hess}

\newcommand{\ds}{\displaystyle}
\newcommand{\wg}{\wedge}
\newcommand{\e}{\mathrm{e}}
\newcommand{\Lip}{\mathrm{Lip}}

\allsectionsfont{\normalsize\sc\center}

\title{\sf 
Space-time Wasserstein controls and 
Bakry-Ledoux type gradient estimates
}
\author{
Kazumasa Kuwada%
\footnote{%
Partially supported 
by the Grant-in-Aid for Young Scientists (B) 22740083}
}
\begin{document}
\maketitle
\raggedbottom 
\begin{abstract}
The duality in Bakry-\'Emery's gradient estimates and Wasserstein controls 
for heat distributions is extended to that in refined estimates 
in a high generality. 
As a result, we find an equivalent condition to 
Bakry-Ledoux's refined gradient estimate involving an upper dimension bound. 
This new condition is described as a $L^2$-Wasserstein control 
for heat distributions at different times. 
The $L^p$-version of those estimates are studied 
on Riemannian manifolds via coupling method. 
\end{abstract}
\begin{list}{\textbf{Key words:}}
{\setlength{\labelwidth}{72pt}
\setlength{\leftmargin}{72pt}}
\item
gradient estimate, Wasserstein distance, heat distribution, 
Hopf-Lax semigroup, Ricci curvature
\end{list}
% \begin{list}{\textbf{Mathematics Subject Classification (2000):}}
% {
%   \setlength{\labelwidth}{72pt}
%   \setlength{\leftmargin}{72pt}
% }
% \item
% \end{list}

% \section{Introduction} 
% \label{sec:intro} 
% 
% Under construction. 

\tableofcontents

\section{Introduction}

Since the pioneering work of Bakry and \'Emery \cite{BE:diff-hyp}, 
their ($L^2$-)gradient estimate 
\begin{equation} \label{eq:BE0}
| \nabla P_t f |^2 \le \e^{-2Kt} P_t ( | \nabla f |^2 ) 
\end{equation}
for a diffusion semigroup $P_t = \e^{t \mathscr{L}}$ 
has been played a prominent role in geometric analysis of 
the diffusion generator $\mathscr{L}$ 
(see \cite{Bak97,Bakry:2006ul,BGL_book,Led_geom-Markov}, for instance). 
Among several further developments, the following refined 
form is studied by Bakry and Ledoux \cite{Bakry:2006vh} 
(see \cite{Wang_equivCD} also): 
For fixed parameters $K \in \R$ and $N \in ( 0 , \infty ]$, 
\begin{equation} \label{eq:BL0}
| \nabla P_t f |^2 \le \e^{-2Kt} P_t ( | \nabla f |^2 ) 
- \frac{ 1 - \e^{-2Kt} }{ NK } ( \mathscr{L} P_t f )^2 . 
\end{equation}
Recently, in \cite{K9,K13}, 
it is revealed that \eqref{eq:BE0} 
is equivalent to the following estimate 
concerning a Lipschitz type bound of heat distributions 
with respect to the $L^2$-Wasserstein distance $W_2$, 
which we call an $L^2$-Wasserstein control, 
in a fairly general situation: 
For $t > 0$ and two probability measures $\mu_0$ and $\mu_1$ 
on the state space, 
\begin{equation} \label{eq:W0}
W_2 ( P_t^* \mu_0 , P_t^* \mu_1 ) \le \e^{-Kt} W_2 ( \mu_0 , \mu_1 ). 
\end{equation}
% It can be regarded as a sort of $L^2/L^2$-duality.
The purpose of this article is to extend such a duality 
by introducing a new inequality like \eqref{eq:W0} 
which corresponds to \eqref{eq:BL0}. 
As we will see in Corollary~\ref{cor:BL<->W2} 
as a special case of our result, 
the estimate \eqref{eq:BL0} is equivalent to the following 
space-time $L^2$-Wasserstein control in an abstract framework: 
For $s, t > 0$ and two probability measures $\mu_0$ and $\mu_1$ 
on the state space, 
\begin{equation} \label{eq:W2}
W_2 ( P_s^* \mu_0 , P_t^* \mu_1 )^2  
\le 
\abra{ 
  \frac{1}{J_{N} ( [ s, t ])} 
  \int_s^t \e^{Kr} J_N (dr) 
}^{-2} 
W_2 ( \mu_0 , \mu_1 )^2 
+ J_N ( [ s , t ] )^2 , 
\end{equation}
where 
\begin{equation} \label{eq:time}
J_N ( A ) 
: = \int_{A} \sqrt{ \frac{NK}{ \e^{2Kr} - 1 } } dr
\end{equation}
for measurable $A \subset [ 0, \infty )$ and hence 
\begin{equation*}
J_N ( [ s, t ] ) 
=  
\begin{cases}
\ds 
\sqrt{\frac{N}{K}} \abra{ \cos^{-1} ( \e^{-Kt} ) - \cos^{-1} ( \e^{-Ks} ) }
& ( K > 0 ), 
\\
\sqrt{2N} ( \sqrt{t} - \sqrt{s} ) 
& ( K = 0 ), 
\\
\ds
\sqrt{\frac{N}{-K}} \abra{ \cosh^{-1} ( \e^{-Kt} ) - \cosh^{-1} ( \e^{-Ks} ) }
& ( K < 0 ).
\end{cases}
\end{equation*} 
Here the function $( \e^{2Kt} - 1 ) / K$ is regarded as $2t$ when $K=0$. 

% Even when $N = \infty$, \eqref{eq:BL0} 
% is a source of several important functional inequalities 
% such as logarithmic Sobolev inequality, 
% Talagrand inequality and Poincar\'e inequality etc. 
% On the other hand, 
In the Bakry-Ledoux gradient estimate \eqref{eq:BL0}, 
the parameters $K$ and $N$ play the role of lower Ricci curvature bound 
and upper dimension bound. 
% is tightly connected with 
% the notion of the combination of lower Ricci curvature bound and upper dimension bound. 
Indeed, under the condition 
\begin{equation} \label{eq:Gamma}
| \nabla f |^2 = \frac12 \mathscr{L} ( f^2 ) - f \mathscr{L} f, 
\end{equation}
the Bakry-Ledoux gradient estimate 
% \eqref{eq:BL0} 
is equivalent (at least formally) to the following inequality, called 
Bakry-\'Emery's curvature-dimension condition or 
Bochner's inequality 
\begin{equation} \label{eq:Gamma2}
\frac12 \mathscr{L} ( | \nabla f |^2 ) 
- \dbra{ \nabla f , \nabla \mathscr{L} f }
\ge 
K | \nabla f |^2 + \frac{1}{N} | \mathscr{L} f |^2 . 
\end{equation}
Note that \eqref{eq:Gamma} is the \emph{definition} of squared norm of gradient 
or carr\'e du champ $| \nabla f |^2 = \dbra{ \nabla f , \nabla f }$
in Bakry-\'Emery theory, 
as in \cite{Bak97,Bakry:2006ul,BGL_book,Led_geom-Markov}. 
On complete Riemannian manifolds with $\mathscr{L} = \Delta$, 
the Bochner-Weitzenb\"ock formula implies that 
Bakry-\'Emery's curvature-dimension condition 
is equivalent to the combination of $\Ric \ge K$ and $\dim \le N$ 
(when $N = \infty$, the latter condition always holds). 
Moreover, even in an abstract framework, 
\eqref{eq:Gamma2} has provided several extensions of 
results in Riemannian geometry concerning these bounds. 
Thus we could say that \eqref{eq:BE0} or \eqref{eq:BL0} 
is placed at the intersection of geometry and analysis. 
While \eqref{eq:BE0} can be applied in a broader situation
such as analysis on infinite dimensional spaces, 
\eqref{eq:BL0} provides qualitatively sharper results 
and hence obtaining \eqref{eq:BL0} or \eqref{eq:Gamma2} 
for $N < \infty$ would be important 
(see \cite{Bak97,Bakry:2006ul,BGL_book,Led_geom-Markov} for instance). 

The Wasserstein control condition \eqref{eq:W0} or \eqref{eq:W2}
serves us a new approach to \eqref{eq:BE0} or \eqref{eq:BL0}
especially on non-smooth spaces. 
% It has been a significant question when 
% \eqref{eq:BL0} or \eqref{eq:Gamma2} is available 
% on non-smooth spaces. 
Since the Bochner-Weitzenb\"ock formula is not available in such a case, 
it was completely unclear when \eqref{eq:Gamma2} holds. 
For this problem, a connection with an alternative formulation of 
``$\Ric \ge K$ and $\dim \le N$'' by optimal transportation 
\cite{Bacher:2010bp,Lott-Vill_AnnMath09,Sturm_Ric,Sturm_Ric2}
has been investigated recently. 
Since those new conditions are stable under geometric operations 
such as the measured Gromov-Hausdorff limit, 
the same stability holds for \eqref{eq:Gamma2} 
once we prove the equivalence between them. 
This equivalence is finally established 
by Ambrosio, Gigli, Savar\'e, Mondino and Rajala 
\cite{Ambrosio:2012tp,AGS_BE-CD,AGS2,AGS3} 
when $N = \infty$ and 
by Erbar, Sturm and 
the author \cite{Erbar:2013wf} when $N < \infty$. 
For connecting Bakry-\`Emery theory based on \eqref{eq:Gamma2} 
with optimal transport approach, the estimate \eqref{eq:W0} 
or \eqref{eq:W2} works as a bridge, 
though what we actually used when $N < \infty$ is \eqref{eq:sWc} below. 
% which is similar to \eqref{eq:W2} but different. 

The emphasis of the result of this paper is put on the fact that 
the equivalence between \eqref{eq:BL0} and \eqref{eq:W2} 
can be extended to more general situation where 
any kind of known curvature-dimension conditions corresponding to 
``$\Ric \ge K$ and $\dim \le N$'' may not hold 
(see Theorem~\ref{th:Duality}). 
For instance, previous results can be applied to 
obtain an estimate like \eqref{eq:W0} 
from an estimate like \eqref{eq:BE0} 
for sub-elliptic diffusions \cite[Section 4]{K9} 
(see \cite[Section 6]{K13} also for other examples). 
We can expect a similar result also 
in the present case as a future application. 
As another kind of generality, we can obtain $L^p/L^{p_*}$-duality, 
as we did in \cite{K9,K13}. 
Actually, an $L^p$-type estimate of \eqref{eq:W2} 
(see \eqref{eq:Wp}) holds on complete Riemannian manifolds, 
and we can obtain the $L^{p_*}$-analog of \eqref{eq:BL0} by our result. 
It is not yet known whether we can obtain the same estimate 
on metric measure spaces satisfying the Riemannian curvature-dimension 
condition in \cite{Erbar:2013wf}, where the $L^2$-estimate holds. 
Note that, in the case $N = \infty$, the $L^p$-version of \eqref{eq:W0} 
follows from \eqref{eq:W0} itself (see \cite{Savare:2013wa}). 
Such a precision has various applications in that case 
and hence it would be interesting to obtain an $L^p$-estimate 
associated with the curvature-dimension condition as well as 
to investigate applications of it. 

The essential idea of the proof of the duality 
in \eqref{eq:BL0} and \eqref{eq:W2} 
is inherited from the earlier studies in \cite{K9,K13}. 
There we regard \eqref{eq:BE0} as a differentiation of \eqref{eq:W0} 
in space variable with a change of viewpoint 
from the space of measures to the space of functions. 
% in duality. 
Then the opposite implication is regarded as an integration 
in space variable, and it is realized by using the Kantorovich duality 
and the analysis of the associated Hopf-Lax semigroup. 
Thus it was natural in those argument that we used no time-dependency 
on the constant and the Markov kernel $P_t$. 
In the present case, the additional term in \eqref{eq:BL0} involving $N$ 
can be regarded as a differentiation in time parameter. 
The reason why the ``integrated'' estimate \eqref{eq:W2} 
deals with distribution of diffusions at different times 
is based on this fact. 
% More precisely, we can say that 
% information of the additional term in \eqref{eq:BL0} 
% is encoded in the difference in time of \eqref{eq:W2}. 
Actually, if we apply \eqref{eq:W2} when $s=t$ 
by taking a limit $s \to t$, it becomes 
the reduced control \eqref{eq:W0}, which is in duality with \eqref{eq:BE0}. 
In the proof of our main theorem, 
we will couple a space parameter with a time parameter. 
As a result, the argument becomes more complicated 
compared with previous ones in \cite{K9,K13}. 
This fact also makes it unclear what is the optimal bound of 
space-time Wasserstein control of type \eqref{eq:W2}. 
Indeed, we obtain \eqref{eq:W2} by choosing a 
(possibly not optimal but admissible) 
space-time reparametrizations in a variational problem arising 
in Proposition~\ref{prop:WC-var} below 
(see Remark~\ref{rem:improve}). 
When we are working with the specified estimate \eqref{eq:BL0}, 
the space-time $W_2$-control \eqref{eq:sWc} involving comparison functions 
seems to be optimal (see Remark~\ref{rem:sWc}). 

$L^p$-type estimates on Riemannian manifolds are obtained 
by stochastic analytic techniques. We construct a variant of 
coupling by parallel transport of diffusion processes 
with different time scales for deriving Wasserstein controls. 
For the construction, we must avoid technical difficulties 
arising from the presence of the cut locus. 
To overcome it, we employ the approximation of the diffusion process 
by geodesic random walks developed in \cite{K8,K10,K12,Renes_poly}. 
The Wasserstein control we will obtain directly 
from our coupling method is slightly weaker than expected 
and we establish an $L^p$-variant of Bakry-\'Emery theory 
(Proposition~\ref{prop:BEp}) in order to derive a sharper result. 
Since we consider a possibly non-symmetric diffusion generator, 
we relies on stochastic analytic techniques again 
to avoid technical difficulties there. 
It might be possible to prove it in a more straightforward way 
without coupling methods. See Remark~\ref{rem:BEp} for an observation. 

As a related result, the implication from \eqref{eq:BL0} to \eqref{eq:W2} 
is obtained by Bakry, Gentil and Ledoux in \cite{BGL} 
by a different method when $K=0$. 
Bolley, Gentil and Guillin \cite{Bolley:2013ts} recently obtained 
a sort of contraction bound of $L^2$-Wasserstein distance 
for two distributions of diffusion at the \emph{same} time 
involving the dimension parameter $N$ from \eqref{eq:Gamma2}. 
Coupling method on Riemannian manifolds is studied first 
by Kendall \cite{Kend} and improved by Cranston \cite{Crans}. 
Since then, it has been extensively studied in the literature. 
Feng-Yu Wang is one of the leading persons on this topic and 
we refer to his books \cite{Wang_book05,Wang_book14} 
for further developments. 
A coupling admitting different time scales is studied also 
in \cite{K-Phili2,NP}. 
The problem studied in \cite{K-Phili2} seems to be 
closely related to ours (See Remark~\ref{rem:RF}). 
Though the purpose of \cite{NP} is different from ours, 
there appears a similar argument. 

The organization of the paper is as follows. 
In the next section, we give a precise definition of our framework 
and state the main theorems.
Since we will deal with the time-evolution of the Markov kernel, 
we discuss it under two different family of assumptions  
in Theorem~\ref{th:Duality} and Theorem~\ref{th:Duality2} 
respectively. 
In those theorems, only a weaker duality is obtained 
and we add a technical assumption 
(strong Feller property of $P_t$) in both cases 
to obtain the full duality result including the relation 
between \eqref{eq:BL0} and \eqref{eq:W2} 
(Corollary~\ref{cor:Duality}). 
The proof of main results except Theorem~\ref{th:BL-Wp} 
is given in Section~\ref{sec:dual}. 
In the proof, we will show two key propositions: 
Proposition~\ref{prop:diff} and Proposition~\ref{prop:WC-var}. 
Though they can be applied to more general situation than 
in the setting of Theorem~\ref{th:Duality}, 
we exclude it from the main theorem 
for simplicity of presentation 
since the statements of them look more complicated. 
We also prove another space-time $W_2$-control 
\eqref{eq:sWc} studied in \cite{Erbar:2013wf} 
directly from \eqref{eq:BL0} (Theorem~\ref{th:sWc}). 
In Section~\ref{sec:coupling}, we will prove 
Theorem~\ref{th:BL-Wp}, which concerns with $L^p$-estimates, 
on a complete Riemannian manifold satisfying \eqref{eq:Gamma2} 
for a (possibly non-symmetric) diffusion generator $\mathscr{L}$. 
% For the proof, we will use a stochastic analytic technique 
% based on a coupling of diffusion processes constructed by approximation. 
Since the argument seems to be technical, 
we give a heuristic discussion in Section~\ref{sec:absence} 
and make it rigorous in Section~\ref{sec:RW} 
with a partial use of arguments in Section~\ref{sec:absence} 
which hold in a sufficient generality. 
\medskip

\noindent
\emph{Acknowledgment}: 
The author would like to tell his gratitude to the anonymous referee. 
His/Her comments help the author to improve the quality of the paper. 
Especially, Proposition~\ref{prop:BEp} and Remark~\ref{rem:BEp} 
are essentially due to the comment. 

\section{Framework and main results} \label{sec:frame}

Let $( M, d )$ be a Polish metric space. 
In this paper we always assume that $d$ is a geodesic metric. 
It means that, for each $x, y \in M$, 
there is a curve $\gm : [ 0 , 1 ] \to M$ 
such that $\gm (0) = x$, $\gm (1) = y$ and 
$d ( \gm (s) , \gm (t) ) = | s - t | d ( x, y )$ 
for $s,t \in [ 0, 1 ]$. We call such a curve $\gm$ 
minimal geodesic joining $x$ and $y$. 
Let $P_t ( x , \cdot ) \in \sP (M)$, $t \ge 0$, $x \in M$  
be a semigroup of Markov kernels on $M$, 
where $\sP (M)$ is the space of all Borel probability measures on $M$. 
We denote the action of the Markov kernel $P_t$ to $f : M \to \R$ 
bounded and measurable by $P_t f$. 
Similarly, the dual action of $P_t$ to $\mu \in \sP (M)$
is denoted by $P_t^* \mu$. 
We denote the set of (bounded) Lipschitz functions 
by $C^{\Lip} (M)$ and $C^{\Lip}_b (M)$ respectively. 
Let us denote the local Lipschitz constant of $f \in C^{\Lip} (M)$ 
at $x$ by $| \nabla f | (x)$. That is, 
\begin{equation*} %\label{eq:local-Lip}
| \nabla f | (x) = \limsup_{y \to x} \frac{| f(y) - f(x) |}{d(x,y)}. 
\end{equation*}
Let $\Lip (f)$ stand for the (global) Lipschitz constant of $f$. 
Recall that we say a curve $( \gm (r) )_{r \in [ 0, 1 ]}$ 
in a metric space $( Y, d_Y )$ is called 
absolutely continuous if there is a nonnegative integrable function 
$\ph$ on $[ 0 , 1 ]$ such that $d_Y ( \gm (s) , \gm (t) ) \le \int_s^t \ph (r) dr$ 
for any $0 \le s \le t \le 1$. 
$\ph$ can be chosen to be the metric derivative 
$| \dot{\gm} |$ given by 
\begin{equation} \label{eq:m-deriv}
| \dot{\gm} | (r) 
: = 
\limsup_{s\to r} 
\frac{ d_Y ( \gm (s) , \gm (r) ) }{ | s - r | }
\end{equation}
(see \cite[Theorem~1.1.2]{AGS}). 
For $1 \le p < \infty$, 
we denote the $L^p$-Wasserstein (pseudo-)distance 
on $\sP (M)$ by $W_p$. 
That is, 
\begin{align*}
W_p ( \mu , \nu ) 
& := 
\inf 
\bbra{ 
  \| d \|_{L^p ( \pi )} 
  \; | \; 
  \mbox{$\pi$ is a coupling of $\mu$ and $\nu$} 
} . 
\end{align*}
For each $\mu_0 , \mu_1 \in \sP (M)$, 
there is a curve $( \mu(t) )_{t \in [ 0 , 1 ]} \subset \sP (M)$ 
such that $\mu (i) = \mu_i$ for $i = 0 , 1$ and 
$W_p ( \mu (r_1) , \mu (r_2) ) = | r_1 - r_2 | W_p ( \mu_0 , \mu_1 )$ 
for any $r_1 , r_2 \in [ 0 , 1 ]$ 
(see \cite[Corollary 1 and Proposition 1]{Lisini:2007be}). 
Let $\mathrm{Geo} (M)$ be the space of minimal geodesics 
parametrized by $[0,1]$ and $e_t : \mathrm{Geo} (M) \to M$ 
($t \in [ 0, 1 ]$) the evaluation map given by 
$e_t (\gm) = \gm (t)$. 
By \cite[Corollary 1 and Theorem 6]{Lisini:2007be}, 
there exists a probability measure $\Gm$ on $\mathrm{Geo} (M)$ 
such that $( e_r )_\sharp \Gm = \mu (r)$ and 
\begin{equation} \label{eq:dynamic}
\int_{\mathrm{Geo}(M)} d ( \gm (r_1) , \gm (r_2) )^p \Gm ( d \gm ) 
= 
\int_{\mathrm{Geo}(M)} \int_{r_1}^{r_2} | \dot{ \gm } | (u)^p du \Gm ( d \gm ) 
= W_p ( \mu (r_1) , \mu (r_2 ) )^p
\end{equation}
for any $r_1 , r_2 \in [ 0 , 1 ]$ with $r_1 < r_2$. 
we call such $\Gm$ a dynamic optimal coupling of $\mu_0$ and $\mu_1$. 

We introduce some quantities we will use throughout this paper. 
Let $a : [ 0 , \infty ) \to ( 0 , \infty )$ and 
$b : ( 0 , \infty ) \to ( 0 , \infty )$ be continuous functions. 
We define a measure $J$ on $[ 0 , \infty )$ 
by $J ( dx ) = b (x)^{-1} dx$. We assume that 
$J$ is locally finite, 
that is, $J ( [ 0, \dl ) ) < \infty$ for any $\dl > 0$. 
Let $p, p_* , \b , \b_* \in ( 1, \infty )$ with $p^{-1} + p_*^{-1} =1$, 
$\b^{-1} + \b_*^{-1} = 1$ and $\b \le p$. 
For $f : M \to \R$, we define the Hopf-Lax (or Hamilton-Jacobi) 
semigroup $( Q_s f )_{s \ge 0}$ by 
\begin{equation} \label{eq:HL}
Q_s f (x) 
: = 
\inf_{y \in X} 
\left[ 
f (y) + \frac{s}{p} \abra{  \frac{d (x,y)}{s} }^p 
\right]. 
\end{equation}

For the infinitesimal generator $\sL$ of $P_t$, 
we suppose either of the following conditions: 
\begin{description}
\item[(A1)] %\label{i:ptwise}
For any $f \in C^{\Lip}_b (M)$, $t > 0$ and $x \in M$, 
the following limit exists: 
\begin{equation} \label{eq:generate}
\sL P_t f (x) : = 
\lim_{s \to 0} \frac{ P_{t+s} f (x) - P_t f (x) }{s} . 
\end{equation}
\item[(A2)] %\label{i:Lp} 
There is a locally finite reference measure $\mathfrak{m}$ on $M$ 
with $\mathop{\mathrm{supp}} \mathfrak{m} = M$ 
such that we can extend the action of $P_t$ to $L^{q} (\mathfrak{m})$ 
as a bounded operator for some $q \in [ 1, \infty )$ 
and the limit \eqref{eq:generate} exists in $L^{q} (\mathfrak{m})$ 
for any $f \in L^{q} (\mathfrak{m})$ and $t > 0$. 
\end{description}
The condition \textbf{(A1)} seems more restrictive, 
but the other assumptions can be rather weak 
and the proof of the main theorem is simpler under this condition. 
The condition \textbf{(A2)} requires some additional assumptions 
for the main theorem, but it naturally occurs 
when we are following a functional analytic approach. 
Such a situation arises in analysis on metric measure spaces 
where no (usual) differentiable structure is assumed. 
Under \textbf{(A2)}, for two measurable functions $f$ and $g$ 
which belong to 
the same equivalence class in $L^q (\mathfrak{m})$, 
$\int_M f(y) d P_t (x,dy) = \int_M g(y) P_t (x,dy)$ 
holds $\mathfrak{m}$-a.e.~$x \in M$. 
Thus, even under \textbf{(A2)}, we always regard $P_t f$ for $f \in L^q (\mathfrak{m})$ 
as the integral by the Markov kernel of a representative of $f$. 

We are interested in the following conditions: 
\begin{description}
\item[(1)] \label{i:W}
(Space-time $( L^p, L^\b )$-Wasserstein control)
For $\mu_0 , \mu_1 \in \sP (M)$ and $0 \le s < t$, 
\begin{align} \label{eq:W}
W_p ( P_s^* \mu_0 , P_t^* \mu_1 )^\b 
\le 
\abra{ \frac{1}{J ( [ s, t ])} \int_{[s,t]} \frac{ J (dr) }{ a (r) } }^{-\b}
W_p ( \mu_0 , \mu_1 )^\b 
+ J ( [ s, t ] )^\b . 
\end{align}

\item[(2)] \label{i:BL}
($(L^{p_*} , L^{\b_*} )$-Bakry-Ledoux type gradient estimate) 
For $f \in C_b^{\Lip} (M)$, $t > 0$ and $x \in M$, 
\begin{align} \label{eq:BL}
| \nabla P_t f | (x)^{\b_*} 
\le 
a (t)^{\b_*} \abra{ 
  P_t ( | \nabla f |^{p_*} ) (x)^{\b_*/p_*}
  -  
  b(t)^{\b_*} | \sL P_t f (x) |^{\b_*} 
}. 
\end{align}
\item[(2)${}^*$] \label{i:Q-BL}
\eqref{eq:BL} holds for 
$t > 0$, $x \in M$ and $f$ of the form $f = Q_\dl \tilde{f}$ 
with $\dl > 0$ and $\tilde{f} \in C_b^{\Lip} (M)$. 
\item[(3)] \label{i:BL-int}
For any minimal geodesic $\gm : [ 0 , 1 ] \to M$, 
$0 \le s \le t$ and $f \in C_b^{\Lip} (M)$, 
\begin{multline} \label{eq:BL-int}%\nonumber
\abs{ P_{t} f ( \gm (1) ) - P_{s} f ( \gm (0) ) } 
\\ 
\le \int_0^1 
\abra{ 
  a( \xi(r) )^\b d(\gm (0), \gm(1))^{\b} 
  + 
  \abra{ 
    \frac{ t-s }{ b (\xi(r)) } 
  }^{\b}  
}^{1/\b}
P_{\xi(r)} 
\abra{ 
    | \nabla f |^{p_*} 
  } 
  ( \gm (r) )^{1/p_*}
dr ,
\end{multline}
where $\xi (r) := rt + (1-r)s$. 
% For any absolutely continuous curve 
% $( \xi (r), \gm (r) )_{r \in [ 0, 1 ]}$ in $( 0 , 1 ) \times M$ 
% and $f \in C_b^{\Lip} (M)$, 
% \begin{multline} \label{eq:BL-int}%\nonumber
% \abs{ P_{\xi(1)} f ( \gm (1) ) - P_{\xi(0)} f ( \gm (0) ) } 
% \\ 
% \le \int_0^1 
% \abra{ 
%   \abra{ a(\xi(r)) |\dot{\gm}| (r) }^{\b} 
%   + 
%   \abra{ 
%     \frac{ |\xi'(r)| }{ b (\xi(r)) } 
%   }^{\b}  
% }^{1/\b}
% P_{\xi (r)} 
% \abra{ 
%     | \nabla f |^{p_*} 
%   } 
%   ( \gm (r) )^{1/p_*}
% dr .
% \end{multline}
\end{description}
Note that \textbf{(2)} implies \textbf{(2)}${}^*$ since the function $f$ in \textbf{(2)}${}^*$ 
belongs to $C_b^{\Lip} (M)$ (see Section~\ref{sec:HL}). 
For \textbf{(2)} and \textbf{(2)}${}^*$, we consider a slightly modified version under \textbf{(A2)}. 
When \textbf{(2)} or \textbf{(2)}${}^*$ holds for 
$f \in C_b^{\Lip} (M) \cap L^q (\mathfrak{m})$ 
and $\mathfrak{m}$-a.e.~$x \in M$ 
instead of $f \in C_b^\Lip (M)$ and $x \in M$, 
we denote those conditions by 
\textbf{(2)}${}_{\mathrm{ae}}$ or \textbf{(2)}${}_{\mathrm{ae}}^*$ respectively. 
Note that, when considering \textbf{(2)}$^*$, 
$f = Q_\dl \tilde{f} \in C_b^{\mathrm{Lip}} (M)$ is automatic 
and $f \in L^q (m)$ holds 
if $\tilde{f} \in L^q (\mathfrak{m})$ and $\tilde{f} \ge 0$, 
or $\mathop{\mathrm{supp}} \tilde{f}$ is compact.  
We state our first main theorems of this paper as follows: 

\begin{Thm} \label{th:Duality}
Assume {\normalfont \textbf{(A1)}}. 
Then the conditions {\normalfont \textbf{(1)}}, 
{\normalfont \textbf{(2)}${}^*$} and 
{\normalfont \textbf{(3)}} are equivalent. 
\end{Thm}
For considering the corresponding assertion under \textbf{(A2)}, 
we introduce the following additional assumption associated with \textbf{(A2)}: 
\begin{description}
\item[(A3)]
With keeping $\mathfrak{m}$ and $q$ introduced in \textbf{(A2)}, 
for any $\mu_0 , \mu_1  \in \sP (M)$ with bounded supports and 
bounded densities with respect to $\mathfrak{m}$, 
there exists a $W_p$-minimal geodesic $( \mu_r )_{r \in [ 0, 1 ]}$ 
such that $\mu_t \ll \mathfrak{m}$ and the density $\rho_t$ 
satisfies $\int_A \rho_t^{q_*} d \mathfrak{m} < \infty$ 
for each $t \in [ 0, 1 ]$ and bounded $A \in \mathcal{B} (M)$, 
where $q_*$ is the H\"older conjugate of $q$. 
\end{description}
\begin{Thm} \label{th:Duality2}
Assume {\normalfont \textbf{(A2)}}. 
Then the implication 
``{\normalfont \textbf{(1)}} $\Rightarrow$ 
{\normalfont \textbf{(3)}} $\Rightarrow$ 
{\normalfont \textbf{(2)}${}^*_{\mathrm{ae}}$}'' holds. 
In addition, the implication 
``{\normalfont \textbf{(2)}${}^*_{\mathrm{ae}}$} $\Rightarrow$
{\normalfont \textbf{(1)}}'' 
also holds true when {\normalfont \textbf{(A3)}} holds. 
\end{Thm}
For stating the full equivalence involving \textbf{(2)} instead of \textbf{(2)}${}^*$, 
we introduce the following assumption on a regularization property of $P_t$.

\begin{description}
\item[(A4)]
$P_t$ is strong Feller, that is, $P_t f \in C_b (M)$ 
for any $t > 0$ and any $f : M \to \R$ bounded and measurable. 
\end{description}

\begin{Cor} \label{cor:Duality}
\begin{enumerate}
\item \label{i:S1}
Assume {\normalfont \textbf{(A1)}} and {\normalfont \textbf{(A4)}}. 
Then {\normalfont \textbf{(3)}} implies 
{\normalfont \textbf{(2)}}. 
In particular, 
{\normalfont \textbf{(1)}}, 
{\normalfont \textbf{(2)}} and 
{\normalfont \textbf{(3)}} are equivalent. 
\item \label{i:S2}
Assume {\normalfont \textbf{(A2)}} and {\normalfont \textbf{(A4)}}. 
Then {\normalfont \textbf{(3)}} implies 
{\normalfont \textbf{(2)}${}_{\mathrm{ae}}$}. 
In particular, {\normalfont \textbf{(1)}}, 
{\normalfont \textbf{(2)}${}_{\mathrm{ae}}$} 
and {\normalfont \textbf{(3)}} are equivalent 
if {\normalfont \textbf{(A3)}} holds additionally. 
\end{enumerate}
\end{Cor}
As a special case of Corollary~\ref{cor:Duality}, 
we obtain the following: 

\begin{Cor} \label{cor:BL<->W2} 
Let $K \in \R$ and $N \in ( 0, \infty )$. 
% Assume that $P_t$ is strong Feller. 
\begin{enumerate}
\item
Assume {\normalfont \textbf{(A1)}} and {\normalfont \textbf{(A4)}}. 
Then the following are equivalent: 
\begin{enumerate}
\item
\eqref{eq:W2} holds for any $\mu_0, \mu_1 \in \sP (M)$ and $0 < s < t$. 
\item
\eqref{eq:BL0} holds for any $f \in C^{\Lip}_b (M)$, $t > 0$ and $x \in M$. 
\end{enumerate}
\item
Assume {\normalfont \textbf{(A2)}}, {\normalfont \textbf{(A3)}} 
and {\normalfont \textbf{(A4)}}. 
Then the following are equivalent: 
\begin{enumerate}
\item
\eqref{eq:W2} holds for any $\mu_0 , \mu_1 \in \sP (M)$ and $0 < s < t$ 
\item
\eqref{eq:BL0} holds for any $f \in C^{\Lip}_b (M) \cap L^q (\mathfrak{m})$, 
$t > 0$ and $\mathfrak{m}$-a.e.~$x \in M$. 
\end{enumerate}
\end{enumerate}
\end{Cor}
Indeed, we obtain Corollary~\ref{cor:BL<->W2} 
from Theorem~\ref{th:Duality} and Theorem~\ref{th:Duality2} 
with $p = \b = 2$, $a(t) = \e^{-Kt}$ and $b(t) = \sqrt{ ( \e^{2Kt} - 1 ) / (NK) }$. 

The reader may think that it seems difficult 
to specify $a$ and $b$ in \eqref{eq:W} 
when we have a bound of $W_p ( P_s^* \mu_0 , P_t^* \mu_1 )$ involving 
$W_p ( \mu_0 , \mu_1 )$, $t$ and $s$. Even in such a case, 
we can find them by passing through our duality argument. 
See Remark~\ref{rem:improve}. 

We next state a equivalence between \eqref{eq:BL0} 
and another Wasserstein control involving comparison functions. 
It is studied in \cite{Erbar:2013wf} in connection with 
the reduced curvature-dimension condition introduced 
in \cite{Bacher:2010bp}. 
Here we give a more direct proof under a slightly different assumptions. 
Let us introduce comparison functions as follows: 
For $\kappa \in \R$, 
we will define functions 
$\mathfrak{s}_\kappa , \mathfrak{c}_\kappa , \mathfrak{t}_\kappa$ 
on $[ 0 , ( \kappa \vee 0 )^{-1/2} \pi ] \cap [0,\infty)$ 
as follows:  
\begin{align*}
\mathfrak{s}_\kappa (u) 
& : = 
\frac{1}{\sqrt{\kappa}} \sin ( \sqrt{\kappa} u ),
&
\mathfrak{c}_\kappa (u) 
& : = 
\cos ( \sqrt{\kappa} u ),
&
\mathfrak{t}_\kappa (u) 
& : = 
\frac{\mathfrak{s}_\kappa (u)}{\mathfrak{c}_\kappa (u)}. 
\end{align*}
When $\kappa > 0$, there is no problem in this definition. 
When $\kappa = 0$, we extend its definition naturally 
to the limit $\kappa \to 0$ of them, as usual. 
Even when $\kappa < 0$, this definition makes sense 
by regarding trigonometric functions as complex functions. 
They take their values in $\R$ even in this case. 

\begin{Thm} \label{th:sWc}
Let $K \in \R$ and $N \in ( 0, \infty )$. 
When $K > 0$, we suppose that $d(x,y) < \pi \sqrt{ (N-1) / K }$ 
holds for any $x,y \in M$. 
\begin{enumerate}
\item
Assume {\normalfont \textbf{(A1)}} and {\normalfont \textbf{(A4)}}. 
Then \eqref{eq:BL0} holds for $f \in C^{\Lip}_b (M)$ and $t > 0$ 
if and only if the following holds: 
For $0 \le s < t$ and $\mu , \nu \in \sP (M)$, 
\begin{multline} \label{eq:sWc}
\mathfrak{s}_{K/N}^2 \abra{ \frac{ W_2( P^*_s \mu , P^*_t \nu ) }{2} }
\le 
\e^{-K(s+t)} \mathfrak{s}_{K/N}^2 \abra{ \frac{ W_2 ( \mu , \nu ) }{2} }
\\
 + 
\frac{N}{2} \frac{ 1 - \e^{-K(s+t)} }{K(s+t)}
( \sqrt{t}-\sqrt{s} )^2 .
\end{multline}
\item
Assume {\normalfont \textbf{(A2)}}, {\normalfont \textbf{(A3)}} 
and {\normalfont \textbf{(A4)}}. 
Then \eqref{eq:BL0} holds for 
$f \in C^{\Lip}_b (M) \cap L^q (\mathfrak{m})$, $t > 0$ 
and $\mathfrak{m}$-a.e.~$x \in M$ 
if and only if \eqref{eq:sWc} holds 
for $0 \le s < t$ and $\mu , \nu \in \sP (M)$. 
\end{enumerate}
\end{Thm}
Note that the condition on the diameter in the last Theorem 
when $K > 0$ holds in typical situations. 
See e.g.~\cite{BL,Bakry:2006vh} 
and \cite[Remark~3.5 and Corollary~3.7]{Erbar:2013wf} 
(cf.~Remark~\ref{rem:BM-BE}). 

Our final main theorem deals with the case when $M$ is 
a $m$-dimensional complete Riemannian manifold. 
We consider the diffusion process $( ( X (t) )_{t \ge 0} , ( \P_x )_{x \in M} )$ 
generated by $\sL = \Delta + Z$, where $Z$ is a smooth vector field. 
For $K \in \R$ and $N \in [ m , \infty )$, 
we say that the Bakry-\'Emery Ricci tensor associated with $\sL$ 
satisfies the $(K,N)$-curvature-dimension bound if the following holds: 
\begin{equation} \label{eq:Z-CD}
\Ric - ( \nabla Z )^\flat - \frac{1}{N-m} Z \otimes Z \ge K, 
\end{equation}
where $( \nabla Z )^\flat$ is a symmetrization of $\nabla Z$ as $(0,2)$-tensor. 
When $N = m$, we interpret \eqref{eq:Z-CD} as $Z = 0$ and $\Ric \ge K$. 
% We \emph{suppose} 
% $
% P_\cdot f \in 
% C^\infty ( ( 0 , \infty ) \times M ) 
% \cap 
% C ( [ 0 , \infty ) \times M )$ 
% for $f \in C_0^\infty (M)$. 
% Though this assumption typically holds 
% by the hypoellipticity of $\Dl$, we state it as an assumption 
% since $\sL$ can be non-symmetric and 
% the proof of smoothness might be more involved. 
Let $P_t (x , \cdot )$ be given by 
distributions of the diffusion process $X (t)$: 
$P_t ( x , \cdot ) = \P_x \circ X(t)^{-1}$. 
In this case, we will obtain the following: 

\begin{Thm} \label{th:BL-Wp}
Assume $p \ge 2$ and 
\eqref{eq:Z-CD} for some $K \in \R$ and $N \in [ m , \infty )$.  
\begin{enumerate}
\item
For $t > s > 0$ and $\mu_0 , \mu_1 \in \sP (M)$, 
\begin{multline} \label{eq:Wp}
W_p ( P_s \mu_0 , P_t \mu_1 )^2  
\le 
\abra{ 
  \frac{1}{J_{N+p-2} ([s,t])} 
  \int_s^t \e^{Kr} J_{N+p-2} (dr) 
}^{-2} 
W_p ( \mu_0 , \mu_1 )^2 
\\
+ J_{N+p-2} ( [ s , t ] )^2 , 
\end{multline}
where $J_\cdot$ is as defined in \eqref{eq:time}. 
\item
For $t > 0$ and $f \in C^\Lip_b (M)$, 
\begin{equation} \label{eq:BLp}
| \nabla P_t f | (x)^2 
\le \e^{-2Kt} P_t ( | \nabla f |^{p_*} )^{2/p_*} 
- \frac{1 - \e^{-2Kt}}{(N+p-2) K} \abs{ \sL P_t f }^2 .
\end{equation}
\end{enumerate}
\end{Thm}
This result can be regarded as 
an $L^p/L^{p_*}$-version of \eqref{eq:BL0} and \eqref{eq:W2}. 
Note that a similar argument implies an estimate of transportation cost 
involving a comparison function (see Theorem~\ref{th:Lp2}). 

The rest of this section consists of a series of remarks concerning 
with Theorem~\ref{th:Duality} and Theorem~\ref{th:Duality2}, 
including a discussion on sufficient conditions of assumptions in these theorems 
(Remark~\ref{rem:cond_Feller} and Remark~\ref{rem:cond_geod}). 

\begin{Rem} \label{rem:p}
For both \eqref{eq:W} and \eqref{eq:BL}, 
the inequality becomes stronger as $p$ increases (or $p_*$ decreases). 
For \eqref{eq:W}, This is based on the fact that 
\eqref{eq:W} holds for any $\mu_0 , \mu_1 \in \mathcal{P} (M)$ 
if it does for any Dirac measures (see Lemma~\ref{lem:Dirac} below). 
Then the problem is reduced to 
an easy application of the H\"older inequality. 
As a by-product of this observation, 
When both $\mu_0$ and $\mu_1$ are Dirac measure, 
\eqref{eq:W} yields the same estimate 
even when $1 \le p < \b$. 
For \eqref{eq:BL}, this is an easy consequence of 
the H\"older inequality. 
Here we do not require the fact $p_* \le \b_*$. 
Note that these arguments do not require the conclusion of
Theorem~\ref{th:Duality} or Theorem~\ref{th:Duality2}. 
Note also that, in the $L^p / L^{p_*}$-estimates 
in Theorem~\ref{th:BL-Wp}, 
the constant corresponding to $b(t)$ 
does depend on $p$. 
Thus it is not clear whether 
the same implication still holds or not. 
\end{Rem}

\begin{Rem} \label{rem:cond_Feller}
The strong Feller property we assumed in Corollary~\ref{cor:Duality} 
holds for the heat semigroup 
associated with the (quadratic) Cheeger energy functional 
on a metric measure space 
with a Riemannian lower Ricci curvature bound 
(see \cite[Theorem~6.1 (iii)]{AGS3} and 
\cite[Theorem~7.1 (iii)]{Ambrosio:2012tp}). 
More generally, 
when the semigroup is associated with a Dirichlet form, 
it is known that \eqref{eq:BE0} or \eqref{eq:W0} 
yields the strong Feller property 
under some regularity assumptions
\cite[Theorem~3.17]{AGS_BE-CD}. 
Since either \eqref{eq:BE0} or \eqref{eq:W0} 
immediately follows from \eqref{eq:BL0} or \eqref{eq:W2}, 
the strong Feller property is closely related with our conditions. 
\end{Rem}

\begin{Rem} \label{rem:cond_geod}
The assumption {\normalfont \textbf{(A3)}} 
is satisfied for any $q \in [ 1 , \infty )$
if $( M, d, \mathfrak{m} )$ enjoys 
the curvature-dimension condition $\textsf{CD} ( K , \infty )$ 
in the sense of \cite{Sturm_Ric} for some $K \in \R$ 
(see \cite{Rajala:2012jm}). 
Note that, even in this framework, 
our semigroup $P_t$ is not necessarily the one 
studied in \cite{Ambrosio:2012tp,AGS2,AGS3}, 
which is associated with the (quadratic) Cheeger energy. 
\end{Rem}

\begin{Rem} \label{rem:RF}
In \cite{K-Phili2,Topp_Lopt}, 
the monotonicity of normalized $\mathcal{L}$-transportation cost 
between two heat distributions on a backward Ricci flow 
is studied. 
It can be written in the following form: 
For $s_1 < s_2$ and $t_1 , t_2 > 0$ with $t_2 / t_1 = s_2 / s_1$, 
\begin{align*}
\mathcal{T}_{L} ( \mu_{s_1 + t_1} , s_1 + t_1 ; \mu_{s_2 + t_2} , s_2 + t_2 ) 
\le 
\mathcal{T}_{L} ( \mu_{s_1} , s_1 ; \mu_{s_2} , s_2 )
+ 
2 m ( \sqrt{ t_2 } - \sqrt{ t_1 } )^2, 
\end{align*}
where $m$ is the dimension of the manifold, 
$\mu_{t}$ is the heat distribution at time $t$, 
\[
\mathcal{T}_{L} ( \mu , s ; \nu , t ) 
: = 
\inf \bbra{ \left. 
  2 ( \sqrt{ t } - \sqrt{s} ) \int L ( x , s ; y , t ) \pi ( d x d y )
  \; \right| \; 
  \ \mbox{$\pi$: coupling of $\mu$ and $\nu$}
} 
\]
and $L$ is Perelman's $\mathcal{L}$-distance. 
It looks very similar to \eqref{eq:W2} with $K=0$. 
% Actually, we just replace the squared $L^2$-Wasserstein distance $W_2^2$
% with the (rescaled) $\mathcal{L}$-transportation cost $\mathcal{T}_{\mathcal{L}}$. 
\end{Rem}

\begin{Rem} \label{rem:Lap_compare}
When $s=0$, $\mu_0 = \dl_y$ and $\mu_1 = \dl_x$, 
the inequality \eqref{eq:sWc} becomes the following form: 
\begin{align*} 
\mathfrak{s}_{K/N}^2 \abra{ 
\frac{ W_2 ( \dl_y , P_t \dl_x ) }{2} 
} 
& = 
\mathfrak{s}_{K/N}^2 \abra{ 
\frac12 
  \abra{ \int_M d (y , z)^2  P_t ( x, dz ) }^{1/2} 
} 
\\
& \le 
\e^{-Kt} \mathfrak{s}_{K/N}^2 \abra{ \frac{ d ( y, x ) }{2} } 
+ \frac{Nt}{2} + o (t) . 
\end{align*}
By applying the H\"older inequality to bound $W_1$ by $W_2$, 
the last inequality formally implies the following estimate
by taking a derivative at $t = 0$: 
\[
( \sL d ( y , \cdot ) ) (x) \le \frac{N}{\mathfrak{t}_{K/N} (d (x,y) )}. 
\]
It corresponds to the Laplacian comparison theorem 
on complete Riemannian manifolds, but slightly weaker 
(it was sharp if we could replace $N$ with $N-1$). 
By the same argument based on \eqref{eq:W2} instead of \eqref{eq:sWc}, 
the sharp estimate follows when $K = 0$. 
By using an estimate in Theorem~\ref{th:Lp2} below (with $p=2$), 
we can recover the sharp estimate for $K \in \R$, 
but it is shown only on a complete Riemannian manifold 
(see e.g.~\cite[Lemma~2.1]{K11} for a more direct proof of 
the Laplacian comparison theorem for $\sL = \Dl + Z$). 
\end{Rem}

\section{Proof of dualities} \label{sec:dual}

Before going into the proof, we review 
known properties of the Hopf-Lax semigroup. 

\subsection{Reminder of the Hopf-Lax semigroup}
\label{sec:HL}

Recall that the Hopf-Lax semigroup is defined as in \eqref{eq:HL}. 
It is immediate from the definition that 
$Q_s f$ is non-increasing in $s$ and 
\begin{equation} \label{eq:Q-bound}
\inf_{y \in M} f (y) \le Q_s f (x) \le f(x) . 
\end{equation}
When $f$ is bounded, we can easily observe 
\begin{equation} \label{eq:Q-bound3}
Q_s f (x) 
= 
\inf 
\bbra{ 
  f (y) + \frac{s}{p} \abra{ \frac{ d (x,y) }{s} }^p 
  \; \left| \; 
 \begin{array}{l}
  y \in M, 
  \\
  \ds d ( x,y ) \le s \abra{ \frac{ p ( \sup f - \inf f ) }{s} }^{1/p}
  \end{array}
\right. 
} 
\end{equation}
(see e.g.~\cite[Proposition~A.3 (1)]{GRS_HJ}). 
Therefore, if $f$ is bounded with bounded support, 
then $Q_s f$ shares the same property. 
In addition, by virtue of \eqref{eq:Q-bound3}, 
$Q_s f \in C^{\Lip}_b (M)$ holds if $f$ is bounded. 
Again by \eqref{eq:Q-bound3}, we have 
\begin{equation} \label{eq:Q-conti}
\lim_{t \to 0} Q_t f (x) \ge \liminf_{y \to x} f (y),  
\end{equation}
where the limit in the left hand side exists since $Q_t f (x)$ is monotone in $t$. 
This estimate together with \eqref{eq:Q-bound} 
yields that 
$\lim_{s \downarrow 0} Q_s f (x) = f(x)$ holds 
for each fixed $x \in M$ 
if $f$ is lower semi-continuous. 
When $f \in C^\Lip_b (M)$, the same argument as in the proof of 
\cite[Theorem~2.1 (iv)]{BEHM_HJ} yields 
\begin{align} \label{eq:Qs-lip}
| Q_s f (x) - Q_s f (y) | 
& \le 
\Lip (f) d ( x, y ),
\\ \label{eq:Qt-lip}
| Q_{s'} f (x) - Q_s f (x) | 
& \le 
\frac{\Lip (f)^{p_*}}{p_*} | s' - s | 
\end{align}
for each $x , y \in M$, $s,s' > 0$. 
As an important property of $Q_s f$, 
it is a solution to a Hamilton-Jacobi equation 
in the following sense: 
\begin{align} \label{eq:HJ}
\frac{\partial^+}{\partial s} Q_s f (x) 
 =  
- \frac{1}{p_*} | \nabla Q_s f |^{p_*} (x)
\end{align}
for any $x \in M$ and $s > 0$ 
(see \cite{AGS_S,GRS_HJ}, 
\cite[Theorem~3.6 and Theorem~3.8]{K13} and references therein). 
Note that we use the property that $M$ is a geodesic space 
to obtain the equality \eqref{eq:HJ} 
while an inequality ``$\le$'' holds without this assumption. 
By \cite[Lemma~3.3 and Proposition~3.4]{K13} (see \cite{AGS2,AGS_S} also) 
and \eqref{eq:HJ}, the function $\ds x \mapsto | \nabla Q_s f | (x)$ 
is upper semi-continuous. 
This fact works as a sort of regularization of $| \nabla f |$. 

\subsection{From Wasserstein control to gradient estimates}
\label{sec:W->BL}
In this subsection, we will give the proof of 
the implication ``\textbf{(1)} $\Rightarrow$ \textbf{(3)} $\Rightarrow$ \textbf{(2)}${}^*$''
in Theorem~\ref{th:Duality} and the corresponding assertions in 
Theorem~\ref{th:Duality2}, Corollary~\ref{cor:Duality} 
and Theorem~\ref{th:sWc}. 
The argument is separated into Proposition~\ref{prop:diff}, 
Proposition~\ref{prop:diff5/3} and Lemma~\ref{lem:diff2}. 
We will show all these implications at the end of this section. 

\begin{Prop} \label{prop:diff}
Let $A : [0,\infty)^2 \to ( 0 , \infty )$ 
and $B : [0,\infty)^2 \to ( 0 , \infty )$ 
continuous functions 
satisfying $A(s,t) = A (t,s)$ and $B(s,t) = B (t,s)$. 
Assume that $B(s,t) = 0$ if and only if $s=t$. 
In addition, we assume $A(t,t) = a(t)$ and 
\[
b(t) = \lim_{s \to t} \frac{|s-t|}{B (s,t)}
\]
for $t \in [ 0 , \infty )$, 
where $a, b$ is as introduced in Section~\ref{sec:frame}. 
Suppose the following inequality holds: 
\begin{align} \label{eq:Wpw}
W_p ( P_s^* \dl_x , P_t^* \dl_y )^\b 
\le 
A (s,t)^\b 
d ( x , y )^\b 
+ B (s,t)^\b 
\end{align}
for any $x , y \in M$ and $0 \le s < t$. 
Then, for any absolutely continuous curve 
$( \xi (r), \gm (r) )_{r \in [ 0, 1 ]}$ in $( 0 , \infty ) \times M$ 
and $f \in C_b^{\Lip} (M)$, 
\begin{multline} \nonumber %\label{eq:BL-int}
\abs{ P_{\xi(1)} f ( \gm (1) ) - P_{\xi(0)} f ( \gm (0) ) } 
\\ 
\le \int_0^1 
\abra{ 
  \abra{ a(\xi(r)) |\dot{\gm}| (r) }^{\b} 
  + 
  \abra{ 
    \frac{ |\xi'(r)| }{ b (\xi(r)) } 
  }^{\b}  
}^{1/\b}
P_{\xi (r)} 
\abra{ 
    | \nabla f |^{p_*} 
  } 
  ( \gm (r) )^{1/p_*}
dr .
\end{multline}
In particular, 
the condition {\normalfont \textbf{(3)}} holds. 
\end{Prop}
Note that neither \textbf{(A1)} nor \textbf{(A2)} is required 
in Proposition~\ref{prop:diff}. 

\begin{Proof}
We first claim that $( t, x ) \mapsto P_t f (x)$ is locally Lipschitz 
for any $f \in C_b^{\Lip} (M)$. 
Indeed, \eqref{eq:Wpw} yields that 
$P_t^* \dl_z$ is locally Lipschitz in $(t,z)$ 
with respect to $W_p$ by the assumption on $A$ and $B$. 
By the H\"older inequality, the same holds for $W_1$. 
Then the claim follows from the Kantorovich-Rubinstein duality 
(see \cite[Remark~6.5]{book_Vil2} for instance). 
Note that the claim implies that 
$(t,z) \mapsto P_t^* \dl_z$ is continuous on $( 0 , \infty ) \times M$ 
with respect to the topology of weak convergence. 
Let $( \xi (s) , \gm (s) )_{s \in [ 0, 1 ]}$ be 
an absolutely continuous curve in $( 0 , \infty ) \times M$. 
Then our claim implies that 
$P_{\xi(s)} f ( \gm (s) )$ is absolutely continuous 
in $s$. Thus it is differentiable a.e. with respect to 
the Lebesgue measure on $[ 0 , 1 ]$. 

Let $s \in [ 0 , 1 )$ 
where $s \mapsto P_{\xi (s)} f ( \gm (s) )$ 
is differentiable and 
take $\ep > 0$ such that $s + \ep \in [ 0 , 1 ]$. 
Let $\pi_s^\ep \in \sP ( M \times M )$ be 
a minimizer of $W_p ( P_{\xi(s)} \dl_{\gm(s)} , P_{\xi (s+\ep)} \dl_{\gm (s+\ep)} )$. 
Then we have 
\begin{equation} \label{eq:couple}
\abs{ 
  P_{\xi(s+\ep)} f (\gm (s+\ep)) - P_{\xi(s)} f (\gm (s)) 
}
\le 
\int_{M \times M} \abs{ f (z) - f (w) } \pi_s^\ep ( dz dw ). 
\end{equation}
Take $r > 0$, which is specified later, and set 
\[
G_r f (z) 
: = 
\sup_{z' ; \; d(z, z') \in (0,r)} 
\frac{ \abs{f (z) - f (z')} }{d (z,z')}. 
\]
Then we have 
\begin{align} \nonumber
\int_{M \times M} \abs{ f (z) - f (w) } \pi_s^\ep ( dz dw )
& = 
\int_{M \times M} \abs{ f (z) - f (w) } 
1_{ \{ d(z,w) \le r \} } \pi_s^\ep ( dz dw )
\\ \nonumber
& \quad 
+ \int_{M \times M} \abs{ f (z) - f (w) } 
1_{ \{ d(z,w) > r \} } \pi_s^\ep ( dz dw )
\\ \nonumber
& \le 
\int_{M \times M} G_r f (z) d( z, w ) \pi_s^\ep (dz dw)
+ 2 \| f \|_\infty \pi_s^\ep ( d > r ) 
\\ \nonumber
& \le 
P_{\xi(s)} ( ( G_r f )^{p_*} )(\gm(s))^{1/p_*} 
W_p ( P_{\xi(s)} \dl_{\gm(s)} , P_{\xi(s+\ep)} \dl_{\gm(s+\ep)} )
\\ \label{eq:difference}
& \hspace{4em} + 
\frac{ 2 \| f \|_\infty }{r^p} 
W_p ( P_{\xi(s)} \dl_{\gm(s)} , P_{\xi(s+\ep)} \dl_{\gm(s+\ep)} )^p , 
\end{align}
where $\| f \|_{\infty} = \sup_{x \in M} |f (x)|$. 
Let us choose $r = r (s,\ep)$ by
\begin{equation*} 
r(s,\ep) 
:= 
W_p ( P_{\xi(s)} \dl_{\gm(s)} , P_{\xi(s+\ep)} \dl_{\gm(s+\ep)} )^{1/(2p_*)} 
.
\end{equation*} 
Applying \eqref{eq:Wpw} to \eqref{eq:difference} 
and substituting it into \eqref{eq:couple}, 
we obtain 
\begin{align} \nonumber 
& \abs{ P_{\xi(s+\ep)} f (\gm(s+\ep)) - P_{\xi(s)} f (\gm(s)) }
\\ \nonumber
& \le 
P_{\xi(s)} ( ( G_r f )^{p_*} )(\gm(s))^{1/p_*} 
\abra{ 
  A (\xi(s),\xi(s+\ep))^\b d(\gm(s),\gm(s+\ep))^\b 
  + 
  B (\xi(s),\xi(s+\ep))^\b }^{1/\b}
\\
& \qquad \label{eq:diff3}
+ 2 \| f \|_{\infty} 
W_p ( P_{\xi(s)} \dl_{\gm(s)} , P_{\xi(s+\ep)} \dl_{\gm(s+\ep)} )^{(p+1)/2} .
\end{align}
Note that $W_p ( P_{\xi(s)} \dl_{\gm(s)} , P_{\xi(s+\ep)} \dl_{\gm(s+\ep)} ) = O (\ep)$ 
as $\ep \downarrow 0$ by \eqref{eq:Wpw}. 
In particular, $r \to 0$ as $\ep \downarrow 0$. 
We divide \eqref{eq:diff3} by $\ep$ and let $\ep \downarrow 0$. 
Since $G_r f \le \Lip (f) < \infty$, 
the dominated convergence theorem yields 
\begin{align*}
\abs{ 
  \frac{\partial}{\partial s} 
  P_{\xi (s)} f ( \gm (s) ) 
}
\le 
P_{\xi (s)} ( | \nabla f |^{p_*} ) ( \gm (s) )^{1/p_*} 
\abra{ 
  \abra{ a (\xi(s)) | \dot{\gm} | (s) }^\b 
  + 
  \abra{ 
    \frac{|\xi'(s)|}{b(\xi(s))}
  }^\b
}^{1/\b}. 
\end{align*}
Thus the assertion holds 
by integrating the last inequality with respect to $s$ 
on $[ 0 , 1 ]$. 
\end{Proof}

\begin{Prop} \label{prop:diff5/3}
Assume {\normalfont \textbf{(3)}}. 
Then $(t,x) \mapsto P_t f(x)$ is Lipschitz on $[ a, b ] \times M$ 
for any $0 < a < b$. In addition, 
\begin{equation*}
| \nabla P_t f | (x)^{\b_*} 
\le 
a (t)^{\b_*} \abra{ 
  P_t ( | \nabla f |^{p_*} ) (x)^{\b_*/p_*}
  -  
  b(t)^{\b_*} \abs{ 
    \frac{\partial}{\partial t} P_t f (x) 
  }^{\b_*} 
} 
\end{equation*}
holds for $f$ of the form $f = Q_\dl \tilde{f}$ 
with $\dl > 0$ and $\tilde{f} \in C_b^{\Lip} (M)$, 
$x \in M$ and $t > 0$ at which $s \mapsto P_s f (x)$ is differentiable. 
Moreover, 
the same conclusion holds for any $f \in C_b^{\Lip} (M)$
under {\normalfont \textbf{(A4)}}. 
\end{Prop}

\begin{Proof}
Let $f \in C_b^{\Lip} (M)$, 
$t,s > 0$ with $t \neq s$ and $x,y \in M$ with $x \neq y$. 
Take a minimal geodesic $\gm$ in $M$ from $x$ to $y$ 
and let $\xi (r) := (t-s)r + s$. 
Then \eqref{eq:BL-int} yields 
\begin{align} \nonumber
| P_{t} & f ( y ) - P_{s} f ( x ) |
\\ \label{eq:diff-int1}
& \le \int_0^1 
\abra{ 
  \abra{ a(\xi(r)) d(x,y) }^{\b} 
  + 
  \abra{ 
    \frac{ |t-s| }{ b (\xi(r)) } 
  }^{\b}  
}^{1/\b}
P_{\xi (r)} 
\abra{ 
    | \nabla f |^{p_*} 
  } 
  ( \gm (r) )^{1/p_*}
dr . 
\end{align}
Because $| \nabla f | \le \Lip (f)$, 
\eqref{eq:diff-int1} implies the claimed Lipschitz continuity. 
Note that $( t, x ) \mapsto P_t^* \dl_x$ is also continuous 
on $( 0 ,  \infty ) \times M$ with respect to 
the topology of weak convergence 
by the Kantorovich-Rubinstein duality. 
Let $\a \in \R \setminus \{ 0 \}$, 
define $\sg_* \in \{ \pm 1 \}$ by 
\[
\sg_* : =
\begin{cases}
1 & 
\mbox{if 
  $\ds 
    \limsup_{y \to x} 
    \frac{ P_t f (y) - P_t f (x) }{ d (y,x) }
    = 
    | \nabla P_t f | (x)
  $,
} 
\\
-1 & \mbox{otherwise} 
\end{cases}
\]
and set $\a_* := \sg_* \a$. 
Take $y \in M$ satisfying $0 < |\a| d(x,y) < t$ 
and set $s := t - \a_* d (x,y) > 0$. 
Then \eqref{eq:diff-int1} yields 
\begin{multline} \label{eq:diff-int2}
\abs{ \frac{P_{t} f ( y ) - P_{s} f ( x ) }{d(x,y)} }
\\
\le 
\int_0^1 
\abra{ 
  1 
  + 
  \abra{ 
    \frac{ |\a| }{ a (\xi(r)) b (\xi(r)) } 
  }^{\b}  
}^{1/\b}
a (\xi(r)) P_{\xi (r)} 
\abra{ 
    | \nabla f |^{p_*} 
  } 
  ( \gm (r) )^{1/p_*}
dr .
\end{multline}

Now we claim that 
$P_{\xi (r)} 
\abra{ 
    | \nabla f |^{p_*} 
  } 
  ( \gm (r) )^{1/p_*}
$ is upper semi-continuous in $r$ 
if either $f = Q_\dl \tilde{f}$ 
for some $\dl > 0$ and $\tilde{f} \in C_b^{\Lip} (M)$, 
or \textbf{(A4)} holds. 
In the former case, 
$| \nabla f |$ is upper semi-continuous 
as reviewed in Section~\ref{sec:HL}. 
Since $Q_{\dl'} ( - | \nabla f |^{p_*} ) \in C_b (M)$ holds 
for $\dl' > 0$, 
\eqref{eq:Q-bound} and \eqref{eq:Q-conti} yield 
\begin{align*}
\limsup_{r' \to r} 
P_{\xi (r')} 
\abra{ 
    | \nabla f |^{p_*} 
  } 
  ( \gm (r') )^{1/p_*}
& \le 
\limsup_{\dl' \downarrow 0} 
\abra{ 
  \lim_{r' \to r} 
  P_{\xi (r')} 
  \abra{ 
    - Q_{\dl'} ( - | \nabla f |^{p_*} )
  } 
  ( \gm (r') )^{1/p_*}
}
\\
& \le 
P_{\xi(r)} ( | \nabla f |^{p_*} ) ( \gm(r) )^{1/p_*} . 
\end{align*}
In the latter case, 
$P_\ep f \in C_b (M)$ for arbitrarily small $\ep > 0$. 
Since $(t,x) \mapsto P_t^* \dl_x$ is continuous, 
$P_{\xi(r)} ( | \nabla f |^{p_*} ) ( \gm (r) )$ is continuous in $r$. 

To conclude \eqref{eq:BL} from \eqref{eq:diff-int2}, 
we consider the left hand side of \eqref{eq:diff-int2}. 
By our choice of $s,t$ and $\sg_*$, we have 
\begin{align*}
\limsup_{y \to x} 
\sg_* \frac{ P_t f (y) - P_s f (x) }{d(x,y)}
& = 
\limsup_{y \to x} 
\abra{ 
  \sg_* \frac{ P_t f (y) - P_t f (x) }{d(x,y)}
  +
  \a \frac{ P_t f (x) - P_s f (x) }{t-s}
}
\\
& = 
| \nabla P_t f | (x) +
\a \frac{\partial}{\partial t} P_t f (x). 
\end{align*}
Thus we obtain 
\begin{align*}
\a \frac{\partial}{\partial t} P_t f (x) 
+ | \nabla P_t f | (x) 
& \le 
\limsup_{y \to x} 
\abs{ 
  \frac{ P_t f (y) - P_s f (x) }{ d (x,y) } 
}
\\
& \le 
\abra{ 
  1 
  + 
  \abra{ 
    \frac{ |\a| }{ a (t) b (t) } 
  }^{\b}  
}^{1/\b}
a (t) P_t 
\abra{ 
    | \nabla f |^{p_*} 
  } 
  ( x )^{1/p_*}. 
\end{align*}
Then the conclusion follows by optimizing over $\a$. 
\end{Proof}
To deal with the case under \textbf{(A2)}, we prepare the following lemma. 
\begin{Lem} \label{lem:diff2}
Assume {\normalfont \textbf{(A2)}} 
and  that $(t,x) \mapsto P_t f (x)$ is Lipschitz 
on $[ a, b ] \times M$ for any $0 < a < b$. 
Then, for each $t > 0$, 
$P_t f (x)$ is differentiable at $t$ $\mathfrak{m}$-a.e.~$x \in M$. 
\end{Lem}
The proof of this lemma goes in a similar way 
as the one for the corresponding assertion 
in the proof of \cite[Theorem~4.4]{Erbar:2013wf}. 
\begin{Proof}
For each $x \in M$, $t \mapsto P_t f (x)$ is differentiable 
for a.e.~$t$ with respect to the Lebesgue measure. 
The Fubini theorem yields that 
the set $\tilde{I} \subset ( 0 , \infty )$ defined by 
\[
\tilde{I} := 
\bbra{
  t \in ( 0 , \infty ) 
  \; | \; 
  \mbox{$P_t f (x)$ is differentiable at $t$ for $\mathfrak{m}$-a.e.~$x \in M$}
}^c 
\] 
is of null Lebesgue measure. 
The proof will be completed once we prove $\tilde{I} = \emptyset$. 
Let $t \in (0 , \infty)$. Then we have $s \in \tilde{I}^c$ with $s < t$. 
Note that \textbf{(A2)} implicitly yields 
that $P_t ( x, A ) = 0$ holds for $\mathfrak{m}$-a.e.~$x \in M$ 
for any measurable $A \subset M$ with $\mathfrak{m} (A) = 0$. 
Since $( t , x ) \mapsto P_t f (x)$ is Lipschitz, 
the dominated convergence theorem implies 
\[
P_{t-s} \abra{ \frac{\partial}{\partial t} P_s f }
= 
P_{t-s} \abra{ \lim_{h \to 0} \frac{ P_{s+h} f - P_s f }{h} }
= 
\frac{\partial}{\partial t} P_t f 
\]
$\mathfrak{m}$-a.e.~and 
hence $P_t f (x)$ is differentiable at $t$ for $\mathfrak{m}$-a.e.~$x \in M$. 
It means $t \notin \tilde{I}$ and hence the assertion holds. 
\end{Proof}

Now we argue the implications 
``\textbf{(1)} $\Rightarrow$ \textbf{(3)} $\Rightarrow$ \textbf{(2)}${}^*$'' 
in Theorem~\ref{th:Duality}, the corresponding implication 
in Theorem~\ref{th:Duality2} and Corollary~\ref{cor:Duality}. 
In all these cases, the implication 
``\textbf{(1)} $\Rightarrow$ \textbf{(3)}'' 
follows immediately from Proposition~\ref{prop:diff} 
by taking $A$ and $B$ as follows: 
for $s < t$, 
% $\b = p$, 
\begin{align*}
A (s,t) := 
\abra{ 
   \frac{1}{J ( [ s , t ] )} 
   \int_{[ s, t ]} \frac{ J (dr) }{ a(r) } 
}^{-1},
\qquad 
B (s,t) := J ( [ s , t ] ) .
\end{align*}
By \textbf{(A1)}, 
the implications ``\textbf{(3)} $\Rightarrow$ \textbf{(2)}${}^*$'' 
in Theorem~\ref{th:Duality} and 
``\textbf{(3)} $\Rightarrow$ \textbf{(2)}'' in Corollary~\ref{cor:Duality}~(i)
is a direct consequence of Proposition~\ref{prop:diff5/3}. 
Under \textbf{(A2)}, 
for $f \in C_b^{\Lip} (M) \cap L^q (\mathfrak{m})$, 
$\frac{\partial}{\partial t} P_t f (x) = \sL P_t f (x)$ 
holds $\mathfrak{m}$-a.e.~$x \in M$ 
if the derivation in the left hand side is defined 
in the classical sense $\mathfrak{m}$-a.e.
Thus Proposition~\ref{prop:diff5/3} and Lemma~\ref{lem:diff2} yields 
the implications ``\textbf{(3)} $\Rightarrow$ \textbf{(2)}${}^*_{\mathrm{ae}}$'' 
in Theorem~\ref{th:Duality2} and 
``\textbf{(3)} $\Rightarrow$ \textbf{(2)}${}_{\mathrm{ae}}$'' 
in Corollary~\ref{cor:Duality}~(ii). 

Finally, we briefly discuss two implications 
``\eqref{eq:sWc} $\Rightarrow$ \eqref{eq:BL}'' 
in Theorem~\ref{th:sWc} (i) and (ii). 
This can be reduced to similar arguments 
because $\mathfrak{s}_{K/N} (u) \sim u$ as $u \to 0$. 

\subsection{From gradient estimate to Wasserstein controls}
\label{sec:BL->W}

For the rest of the proof of 
Theorem~\ref{th:Duality}, Theorem~\ref{th:Duality2} 
and Theorem~\ref{th:sWc},  
the estimate of the Wasserstein distance 
between Markov kernels given in Proposition~\ref{prop:WC-var} 
and Proposition~\ref{prop:WC-var2} below is essential. 
We begin with the following auxiliary lemma. 
\begin{Lem} \label{lem:Dirac}
Let $0 < t' \le t''$ and $C_1 , C_2 \ge 0$. 
If
\[
W_p ( P_{t'}^* \mu_0 , P_{t''}^* \mu_1 )^\b 
\le C_1 W_p ( \mu_0 , \mu_1 )^\b + C_2
\]
holds whenever $\mu_0 , \mu_1$ are Dirac measures, 
then the same holds for any $\mu_0 , \mu_1 \in \sP (M)$. 
\end{Lem}

\begin{Proof}
The proof goes along the same line as \cite[Lemma~3.3]{K9} 
(cf.~\cite[Theorem~4.8]{book_Vil2}). 
Thus we omit some technical details. 
For each $x_1 , x_2 \in M$, 
take an optimal coupling $\tilde{\pi}_{x_1 x_2} \in \sP (M)$ 
of $P_{t'}^* \dl_{x_1}$ and $P_{t''}^* \dl_{x_2}$. 
Let $\pi \in \sP (M^2)$ be an optimal coupling of $\mu_0$ and $\mu_1$ 
and define $\tilde{\pi} \in \sP (M^2)$ by 
\[
\tilde{\pi} (A) := 
\int_{M^2} \tilde{\pi}_{x_1 x_2} ( A )
\pi ( d x_1  d x_2 ). 
\] 
Then the assumption and 
the Minkowski inequality for $L^{p/\b}$-norm yield 
\begin{align*}
W_p ( P_{t'}^* \mu_0 , P_{t''}^* \mu_1 )^\b 
& \le 
\abra{ 
  \int_{M^2} W_p ( P_{t'}^* \dl_{x_1} , P_{t''}^* \dl_{x_2} )^p \pi ( d x_1 d x_2 )
}^{\b/p}
\\
& \le 
\abra{ 
  \int_{M^2} ( C_1 d ( x_1 , x_2 )^\b + C_2 )^{p/\b} \pi ( d x_1 d x_2 )
}^{\b/p}  
\\
& \le 
C_1 \| d^\b \|_{L^{p/\b} (\pi)} + C_2 
\\
& = C_1 W_p ( \mu_0 , \mu_1 )^\b + C_2. 
\end{align*}
Hence the conclusion holds. 
\end{Proof}

\begin{Prop} \label{prop:WC-var}
Let $0 \le s \le t$. 
Let $\eta : [ 0, 1 ] \to [ 0 , 1 ]$ and  
$\xi : [ 0 , 1 ] \to [ s, t ]$ 
be $C^1$-increasing surjections. 
Assume {\normalfont \textbf{(A1)}} and 
{\normalfont $\mathbf{\textbf{(2)}{}^*}$}. 
Then, for $\mu_0 , \mu_1 \in \sP (M)$, 
\begin{equation} \label{eq:WC-var}
W_p ( P^*_{s} \mu_0 , P^*_{t} \mu_1 )^\b 
\le 
\int_0^1
\left(
  a ( \xi (r) )^\b W_p ( \mu_0 , \mu_1 )^\b \eta' (r)^\b
  + 
  \abra{ 
    \frac{ \xi' (r) }{ b ( \xi (r) ) } 
  }^\b 
\right)
d r. 
\end{equation}
\end{Prop}

\begin{Proof}
By virtue of Lemma~\ref{lem:Dirac}, 
it suffices to show the assertion when 
$\mu_0 = \dl_{x_0} , \mu_1 = \dl_{x_1}$, 
$x_0 ,x_1 \in M$. 
However, for later use, 
we argue with general $\mu_0 , \mu_1 \in \sP (M)$ 
for a while. 

Take a $W_p$-minimal geodesic $( \mu (r) )_{r \in [ 0 , 1 ]}$ from $\mu_0$ to $\mu_1$. 
Let $t' , t'' \in [ 0 , 1 ]$ with $t'' > t'$ and set $h = t'' - t'$. 
Then the Kantorovich duality yields
\begin{align} \nonumber
\frac1p & 
\abra{ \frac{ W_p ( P^*_{\xi(t')} \mu (\eta(t')) , P^*_{\xi(t'')} \mu (\eta(t'')) )}{h} }^p 
\\ \nonumber
& = 
\frac{1}{h} 
\sup_{f \in C^{\Lip}_b (M)} 
\left[
  \int_M Q_h f \, d P^*_{\xi(t'')} \mu (\eta(t'')) 
  - 
  \int_M f \, d P^*_{\xi(t')} \mu (\eta(t')) 
\right]
\\ \label{eq:Kant}
& = 
\frac{1}{h}
\sup_{f \in C^{\Lip}_b (M)} 
\left[
  \int_M P_{\xi(t'')} ( Q_h f ) \, d \mu (\eta(t'')) 
  - 
  \int_M P_{\xi(t')} f \, d \mu (\eta(t'))
\right]. 
\end{align}
Let $\Gm$ be a dynamic optimal coupling 
associated with $( \mu(r) )_{r \in [ 0 , 1 ]}$. 
Note that \eqref{eq:BL} is available for $Q_r f$ instead of $f$ for $r > 0$ 
since we assume \textbf{(2)}${}^*$. 
With keeping this fact in mind, 
for $r_i \in [ 0 , h ]$ and $s_i , t_i \in [ t', t'' ]$ ($i=1,2$), 
we have 
\begin{align} \nonumber
\Bigg|
  & 
  \int_M P_{\xi(t_2)} Q_{r_2} f \, d \mu ( \eta (s_2) )  
  - 
  \int_M P_{\xi(t_1)} Q_{r_1} f \,d \mu ( \eta (s_1) )
\Bigg|
\\ \nonumber
& \quad \le 
\left| 
    \int_{\mathrm{Geo} (M)} 
    \int_{\eta(s_1)}^{\eta(s_2)} 
    |\nabla P_{\xi (t_2)} Q_{r_2} f | ( \gm (s) ) |\dot{\gm} | (s)
    d s 
    \Gm ( d \gm )
\right| 
\\ \nonumber
& \hspace{12em} + 
\left| 
    \int_M 
    \left( 
        P_{\xi(t_2)} Q_{r_2} f  
        - 
        P_{\xi(t_1)} Q_{r_1} f  
    \right) d \mu (\eta(s_1)) 
\right|
\\ \nonumber 
& \quad \le 
\Lip (f) a ( \xi (t_2) ) W_p ( \mu ( \eta(s_2) ) , \mu ( \eta(s_1) ) )
+ 
\frac{\Lip(f)^{p_*}}{p_*} | r_2 - r_1 | 
\\ \nonumber 
& \hspace{18em}
+ 
\left| 
    \int_M \int_{\xi (t_1)}^{\xi (t_2)} 
      \sL P_t Q_{r_1} f 
    \, dt 
    d \mu ( \eta(s_1) ) 
\right|
\\ \nonumber
& \quad \le 
\Lip (f) a ( \xi (t_2) ) W_p ( \mu_0 , \mu_1 ) | \eta(s_2) - \eta(s_1) | 
\\ \label{eq:Lip^3} 
& \hspace{16em}
+ 
\frac{\Lip(f)^{p_*}}{p_*} | r_2 - r_1 | 
+ 
\Lip(f) 
\int_{\xi (t_1)}^{\xi (t_2)} \frac{dt}{b(t)}. 
\end{align} 
Here we used the fact that local Lipschitz constant is 
an upper gradient in the first inequality. 
The second inequality follows from 
\eqref{eq:BL}, 
\eqref{eq:Qs-lip} and \eqref{eq:Qt-lip}. 
The third inequality follows from 
\eqref{eq:BL} and \eqref{eq:Qs-lip} again. 
Therefore 
$(r,s,t) \mapsto \int_M P_{\xi (t)} Q_{r} f \, d \mu (\eta(s))$ 
is continuous on $[0 , h ] \times [ t', t'' ]^2$ and 
locally Lipschitz on $[ 0 , h ] \times [t',t''] \times ( t' , t'' ]$. 
In particular, $r \mapsto \int_M P_{\xi (r+t')} Q_{r} f \, d \mu (\eta(r+t')) $ is 
continuous on $[0,h]$ and 
locally Lipschitz on $( 0 , h ]$. 
Hence 
%, for $t_1 , t_2 \in [ 0 , 1 ]$ with $t_1 < t_2$, 
we can apply \cite[Lemma~4.3.4]{AGS} twice to obtain 
\begin{align} \nonumber
\int_M & P_{\xi (t'')} Q_{h} f \, d \mu (\eta(t'')) 
 - 
\int_M P_{ \xi(t')} f \, d \mu (\eta(t')) 
 = 
\int_{t'}^{t''} \frac{\partial}{\partial r}  
\left( 
    \int_M P_{\xi (r)} Q_{r - t'} f \, d \mu (\eta(r))
\right) 
d r  
\\ \nonumber 
& \le 
\int_{t'}^{t''} \int_{\mathrm{Geo} (M)}
\Bigg( 
  | \nabla P_{\xi (r)} Q_{r - t'} f | ( \gm  ( \eta (r) ) )
  | \dot{\gm} | (\eta(r)) \eta' (r) 
  + 
  \xi' (r) \sL P_{\xi (r)} Q_{r - t'} f ( \gm ( \eta (r) ) )
\\ \label{eq:interpolate1}
& \hspace{14em} 
  - 
  \frac{1}{p_*} 
  P_{\xi (r)} ( | \nabla Q_{r - t'} f |^{p_*} ) ( \gm ( \eta (r) ) )
\Bigg) 
\Gm ( d \gm ) d r  
\end{align}
with the aid of \eqref{eq:HJ}. 
By the H\"older inequality and \eqref{eq:BL}, 
\begin{align*}
| \nabla & P_{\xi (r)} Q_{r - t'} f | ( \gm (r) ) 
| \dot{\gm} | ( \eta (r) ) \eta' (r) 
+ 
\xi' (r) \sL P_{\xi (r)} Q_{r - t'} f ( \gm (r) )
\\
& \le 
\abra{ 
  d ( \gm (0) , \gm (1) )^\b \eta' (r)^\b 
  + 
  \abra{ \frac{\xi' (r)}{a ( \xi (r) ) b ( \xi (r) )} }^{\b}
}^{1/\b}
\\
& \hspace{2em} \times
\abra{
    | \nabla P_{\xi (r)} Q_{r - t'} f | ( \gm (r) )^{\b_*} 
    + 
    a ( \xi (r) )^{\b_*} b ( \xi (r) )^{\b_*}
    \abs{ \sL P_{\xi (r)} Q_{r - t'} f ( \gm (r) ) }^{\b_*} 
}^{1/\b_*}
\\
& \le 
\abra{ 
  a ( \xi (r) )^\b d ( \gm (0) , \gm (1) )^\b \eta' (r)^\b 
  + 
  \abra{ \frac{\xi' (r)}{b ( \xi (r) )} }^\b
}^{1/\b}
P_{\xi (r)} ( | \nabla Q_{r - t'} f |^{p_*} ) ( \gm (r) )^{1/p_*}. 
\end{align*}
By combining the last inequality with \eqref{eq:interpolate1} 
and \eqref{eq:Kant}, we obtain 
\begin{multline} \label{eq:pre_WC-var}
\frac1p 
\abra{ 
  \frac{ W_p ( P^*_{\xi(t')} \mu (\eta(t')) , P^*_{\xi (t'')} \mu (\eta(t'')) ) }{h} 
}^p 
\\
\le 
\frac{1}{ph} \int_{t'}^{t''} \int_{\mathrm{Geo}(M)} 
\left(
  a ( \xi (r) )^\b d ( \gm (0) , \gm (1) )^\b \eta' (r)^\b 
  + 
  \abra{ 
    \frac{ \xi' (r) }{ b ( \xi (r) ) } 
  }^\b 
\right)^{p/\b} 
\Gm ( d \gm ) d r
\end{multline}
by the Hausdorff-Young inequality. 
This estimate yields that 
$r \mapsto P^*_{\xi(r)} \mu (\eta(r))$ is an absolutely continuous curve in $W_p$ and 
its metric derivative, denoted by $| \partial_r ( P^*_{\xi(r)} \mu ( \eta (r) ) ) |_{W_p}$, 
satisfies 
\begin{align} \label{eq:speed0}
| \partial_r ( P^*_{\xi(r)} ( \eta (r) ) ) |_{W_p}^p
\le 
\int_{\mathrm{Geo}(M)} 
  \left(
      a ( \xi (r) )^\b d ( \gm (0) , \gm (1) )^\b \eta' (r)^\b 
      + 
      \abra{ 
        \frac{ \xi' (r) }{ b ( \xi (r) ) } 
      }^\b 
  \right)^{p/\b} 
\Gm ( d \gm )
\end{align}
for almost every $r \in [ 0 , 1 ]$. 
Then the conclusion follows if both $\mu_0$ and $\mu_1$ are Dirac measures 
by using the property of metric derivative reviewed in Section \ref{sec:frame}. 
Indeed, in this case $\Gamma$ is a Dirac measure on a minimal geodesic joining 
points where Dirac masses are taking place and 
thus $p$-th powers in the both side of \eqref{eq:speed0} can be removed. 
As remarked at the beginning of the proof, it is sufficient to complete the proof. 
\end{Proof}

\begin{Prop} \label{prop:WC-var2}
Assume {\normalfont \textbf{(A2)}}, {\normalfont \textbf{(A3)}} 
and {\normalfont \textbf{(2)}${}^*_{\mathrm{ae}}$}. 
Let $\eta$ and $\xi$ be as in Proposition~\ref{prop:WC-var}. 
Then, \eqref{eq:WC-var} holds for $\mu_0 , \mu_1 \in \sP (M)$. 
\end{Prop}

\begin{Proof} 
We closely follow the proof of Proposition~\ref{prop:WC-var}. 
We first show \eqref{eq:pre_WC-var} 
when $\mu_0 , \mu_1 \in \sP (M)$ have bounded supports 
and bounded densities with respect to $\mathfrak{m}$. 
Let $\Gm$ be a dynamic optimal coupling associated 
with a geodesics $( \mu (r) )_{r \in [ 0 , 1 ]}$ 
from $\mu_0$ to $\mu_1$ given as in \textbf{(A3)}. 
In the Kantorovich duality \eqref{eq:Kant}, 
we may restrict the class of test functions $f$ 
to be $C^{\Lip} (M)$ with bounded supports. 
For such $f$, $Q_r f$ belongs to the same class again 
for any $r > 0$. 
In particular, the local finiteness of $\mathfrak{m}$ implies 
$Q_r f \in L^q ( \mathfrak{m} )$. 
Thus the combination of \textbf{(A2)} and 
the choice of $( \mu (r) )_{r \in [ 0 , 1 ]}$ 
make the computation in \eqref{eq:Lip^3} and \eqref{eq:interpolate1} valid.  
Indeed, though \eqref{eq:BL} holds only $\mathfrak{m}$-a.e., 
it is sufficient in this case since $\mu (r) \ll \mathfrak{m}$ by \textbf{(A3)}. 
Then the rest of the proof of Proposition~\ref{prop:WC-var} works 
in exactly the same way. 

Next we take an approximation of Dirac measures. 
By applying \eqref{eq:pre_WC-var} with $t' = t'' = t$ and $\eta (r) = r$, 
we obtain 
\[
W_p ( P_{t}^* \mu_0 , P_{t}^* \mu_1 ) 
\le 
a(t) W_p ( \mu_0 , \mu_1 ). 
\]
By virtue of this estimate, 
\eqref{eq:WC-var} for Dirac measures follows 
by tending $\mu_i \to \dl_{x_i}$ with respect to $W_p$ 
in \eqref{eq:pre_WC-var} for each $x_0 , x_1 \in \sP (M)$. 
\end{Proof}

Now we will show the implication ``\textbf{(2)}${}^*$ $\Rightarrow$ \textbf{(1)}''
in Theorem~\ref{th:Duality} and 
the corresponding implication 
``\textbf{(2)}${}^*_{\mathrm{ae}}$ $\Rightarrow$ \textbf{(1)}'' 
in Theorem~\ref{th:Duality2}. 
We give the proof only to the former one 
since the other proof goes in exactly the same way. 
We apply Proposition~\ref{prop:WC-var} 
with specified $\xi$ and $\eta$. 
Let us define $\Xi : [ s , t ] \to [ 0 , \infty )$ by 
$\Xi (r) : = J ( [ s , r ] )$, where $J$ is given in 
the statement of Theorem~\ref{th:Duality}. 
By using it, we choose $\xi$ and $\eta$ as follows: 
\begin{align*}
\xi (r) 
& : = 
\Xi^{-1} ( ( \Xi (t) - \Xi (s) ) r + \Xi (s) ) , 
\\
\eta (r) 
& : = 
\abra{ \int_0^1 \frac{du}{ a (\xi (u) ) } }^{-1} 
\int_0^r \frac{du} { a (\xi (u)) }. 
\end{align*} 
We can easily verify that $\xi$ and $\eta$ 
satisfy all conditions we supposed 
in Proposition~\ref{prop:WC-var}. 
Moreover, we have 
\begin{align*}
\frac{ \xi' (r) }{b (\xi (r) )} 
& = 
( \Xi \circ \xi )' (r)
\equiv 
\Xi (t) - \Xi (s) 
= 
J ( [ s , t ] ) , 
\\
a ( \xi (r) ) \eta' (r) 
& \equiv 
\abra{ \int_0^1 \frac{du}{ a (\xi (u) ) } }^{-1} 
= 
\abra{ 
  \frac{1}{ J ( [ s, t ] ) } 
  \int_0^1 \frac{\xi' (u) du}{ a (\xi (u) ) b ( \xi (u) ) }
}^{-1} 
\\
& = 
\abra{ 
  \frac{1}{J ( [ s, t ] )} \int_s^t \frac{J (dr)}{a (r)} 
}^{-1}. 
\end{align*}
By substituting them into \eqref{eq:WC-var}, 
we obtain the desired inequality \eqref{eq:W}. 

By putting our arguments together, 
we complete the proof of Theorem~\ref{th:Duality}, 
Theorem~\ref{th:Duality2} and Corollary~\ref{cor:Duality}. 

\begin{Rem} \label{rem:improve}
The combination of 
Proposition~\ref{prop:diff} 
and Theorem~\ref{th:Duality} 
(or Theorem~\ref{th:Duality2}) 
implies that 
the inequality \eqref{eq:Wpw} has 
a self-improvement property. 
That is, starting from a less sharp estimate of 
type \eqref{eq:Wpw}, we can obtain a sharper estimate of 
type \eqref{eq:W} by passing through \eqref{eq:BL}. 

On the other hand, we can easily obtain 
many weaker inequalities of type \eqref{eq:Wpw} 
from \eqref{eq:BL}. Indeed, since our proof is based on 
Proposition~\ref{prop:WC-var} and an appropriate choice of 
space-time reparametrization $\eta$ and $\xi$, 
a bad choice of $\eta$ and $\xi$ produces a weaker estimate. 
Nevertheless, Proposition~\ref{prop:diff} ensures that 
such a weaker estimate reproduces \eqref{eq:BL} 
and consequently such a \emph{weaker} estimate of type \eqref{eq:Wpw} 
\emph{can be equivalent} to \eqref{eq:BL}. 
Indeed, our choice of $\xi$ and $\eta$ in the proof of 
Theorem~\ref{th:Duality} may not be optimal. 
$\xi$ is a minimizer of the right hand side of \eqref{eq:WC-var}
when $\eta' \equiv 0$ and 
$\eta$ is a minimizer of the same quantity for fixed $\xi$. 
Even in the case of Corollary~\ref{cor:BL<->W2}, 
the genuine minimizer $( \xi , \eta )$ seems to be rather complicated 
(cf.~Remark~\ref{rem:sWc}). 
\end{Rem}

By using Proposition \ref{prop:WC-var},
or Proposition~\ref{prop:WC-var2}, 
we also conclude Theorem~\ref{th:sWc}. 

\begin{tProof}{Theorem~\ref{th:sWc}}
We only give the proof of the implication 
``\eqref{eq:BL0} $\Rightarrow$ \eqref{eq:sWc}'' in (i)
since the other implication is already shown at the end of 
Section~\ref{sec:W->BL} and the corresponding assertion in (ii) 
can be shown in the same way by using Proposition~\ref{prop:WC-var2} 
instead of Proposition~\ref{prop:WC-var}. 
To complete the proof it is sufficient to obtain 
the following differential inequality: 
For $u > 0$ and $\lm \ge 1$, 
\begin{multline} \label{eq:W-var}
\frac{\partial \;}{\partial u}
\mathfrak{s}_{K/N} 
\abra{ 
  \frac{W_2 ( P_{\lm^{-1} u}^* \mu , P_{\lm u}^* \nu )}{2}
}^2
\\
\le 
- K ( \lm + \lm^{-1} ) 
\mathfrak{s}_{K/N} 
\abra{ 
    \frac{ W_2 ( P_{\lm^{-1} u}^* \mu , P_{\lm u}^* \nu) }{2}
}^2
+ 
\frac{N}{2} ( \lm + \lm^{-1} - 2 ). 
\end{multline}
Indeed, letting $\lm = \sqrt{t/s}$ and $r = \sqrt{ts}$ and 
applying the Gronwall lemma to give an estimate 
of $\mathfrak{s}_{K/N} ( W_2 ( P_{\lm^{-1} r}^* \mu , P_{\lm r}^* \nu )/2 )^2$ 
yield the conclusion. 
In the sequel, we use the abbreviation
$w := W_2 ( P^*_{\lm^{-1} u} \mu , P^*_{\lm u} \nu )$ 
for simplicity of notation. 
For $h > 0$, let us define $l$, $\te_h$ and $\xi_h$ as follows: 
\begin{align*}
l (r) 
& := 
\begin{cases}
\mathfrak{c}_{K/N}^{-1} ( \e^{-Kr} ) 
& ( K \neq 0 ), 
\\
\sqrt{ 2 N r } 
& ( K = 0 ), 
\end{cases}
\\
\te_h (r) 
& := 
\begin{cases}
\displaystyle
\frac{ 
  l ( \lm h ) \mathfrak{s}_{K/N} ( w r ) 
  + 
  l ( \lm^{-1} h ) \mathfrak{s}_{K/N} ( w (1 - r ) ) 
}
{ 
  \mathfrak{s}_{K/N} ( w )
}
& ( K \neq 0, w \neq 0 ), 
\\
l ( \lm h ) r + l ( \lm^{-1} h ) ( 1 - r )
& (\mbox{$K = 0$ or $w = 0$}), 
\end{cases}
\\
\xi_h (r) 
& : = 
l^{-1} ( \te_h (r) ) 
= 
\begin{cases}
\displaystyle
- \frac{1}{K} \log ( \mathfrak{c}_{K/N} ( \te_h (r) ) ) 
& ( K \neq 0 ), 
\\
\displaystyle
\frac{\te_h (r)^2}{2N}
& ( K = 0 ). 
\end{cases}
\end{align*}
In what follows, we only consider the case $K \neq 0$ and 
$w \neq 0$ for simplicity of presentation. 
Indeed, the same argument also works in other cases. 
Note that $\xi_h$ is $C^1$-increasing surjection 
from $[0,1]$ to $[ \lm^{-1} h , \lm h ]$. 
Thus applying Proposition~\ref{prop:WC-var} with 
$a (t) = \e^{-Kt}$, 
$b (r) = \sqrt{(\e^{2Kt} - 1 )/(NK)}$, 
$t' = \lm^{-1} h$, 
$t'' = \lm h$, 
$\mu_0 =   P^*_{\lm^{-1} u} \mu$, 
$\mu_1 =   P^*_{\lm u} \nu$, 
$\eta (r) = r$ and  
$\xi = \xi_h$, 
we obtain 
\begin{equation} \label{eq:WC-var1}
W_2 ( P^*_{\lm^{-1} ( u + h ) } \mu , P^*_{\lm ( u + h )} \nu )^2 
\le 
\int_0^1
\left(
  \e^{-2K \xi_h (r)} w^2
  + 
  \frac{ N K }{ \e^{2K\xi_h (r)} - 1 } 
  \xi_h' (r)^2 
\right) 
d r. 
\end{equation}
Since we have 
\begin{align*}
\xi_h' (r) & = 
\frac{\mathfrak{t}_{K/N} ( \te_h (r) )}{N} \te_h' (r)
= b ( \xi_h (r) ) \te_h' (r) , 
& 
\lim_{h \to 0} \frac{l ( \a h )}{\sqrt{h}} 
& = 
\sqrt{ 2 N \a }
\end{align*}
for $\a \ge 0$ and 
the addition formulae 
\begin{align*}
\mathfrak{c}_{K/N} ( u + v ) 
& = 
\mathfrak{c}_{K/N} (u) \mathfrak{c}_{K/N} (v)- \frac{K}{N} \mathfrak{s}_{K/N} (u) \mathfrak{s}_{K/N} (v), 
\\
\mathfrak{s}_{K/N} ( u + v ) 
& = 
\mathfrak{s}_{K/N} (u) \mathfrak{c}_{K/N} (v) + \mathfrak{c}_{K/N} (u) \mathfrak{s}_{K/N} (v), 
\end{align*} 
\eqref{eq:WC-var1} implies 

\begin{align} \nonumber
\frac{\partial \;}{\partial u} &
W_2 ( P^*_{\lm^{-1} u} \mu , P^*_{\lm u} \nu )^2 
% \\ \nonumber
% & 
\le 
\limsup_{h \to 0} 
\frac{1}{h} 
\int_0^1
\left(
  ( \e^{-2K \xi_h (r)} - 1 ) w^2
  + 
  \te_h' (r)^2
\right) 
d r
\\ \nonumber
& = 
\int_0^1 
\Bigg\{ 
- 2K w^2 
\left( 
    \frac{ 
      \sqrt{\lm} \mathfrak{s}_{K/N} ( w r ) 
      + 
      \sqrt{\lm^{-1} } \mathfrak{s}_{K/N} ( w ( 1 - r ) ) 
    }
    { \mathfrak{s}_{K/N} ( w ) }
\right)^2 
\\ \nonumber
& \hspace{6em} + 2N w^2 
\left( 
    \frac{ 
      \sqrt{\lm} \mathfrak{c}_{K/N} ( w r ) 
      - 
      \sqrt{\lm^{-1} } \mathfrak{c}_{K/N} ( w ( 1 - r ) ) 
    }
    { \mathfrak{s}_{K/N} ( w ) }
\right)^2 
\Bigg\}
d r 
\\ \nonumber
& =
\frac{ 2N w^2 }{ \mathfrak{s}_{K/N} ( w )^2 }
\left\{ 
( \lm + \lm^{-1} ) \int_0^1 \mathfrak{c}_{K/N} (2wr) d r 
- 2 \int_0^1 \mathfrak{c}_{K/N} ( w ( 2 r - 1 )  ) d r 
\right\}
\\ \label{eq:W-var2}
& = 
\frac{ 4 w }{ \mathfrak{s}_{K/N} (w) } 
\left( 
    - K ( \lm + \lm^{-1} ) \mathfrak{s}_{K/N} \left( \frac{w}{2} \right)^2 
    + \frac{N}{2} ( \lm + \lm^{-1} - 2 ) 
\right). 
\end{align}
Since 
\begin{equation*}
\frac{\partial \;}{\partial u}
\mathfrak{s}_{K/N} \left( 
    \frac12 
    W_2 ( P^*_{\lm^{-1} u} \mu , P^*_{\lm u} \nu ) 
\right)^2 
= 
\frac{\mathfrak{s}_{K/N} (w)}{4w}
\frac{\partial \;}{\partial u} 
W_2 ( P^*_{\lm^{-1} u} \mu , P^*_{\lm u} \nu )^2 ,
\end{equation*}
\eqref{eq:W-var2} immediately yields \eqref{eq:W-var}.  
\end{tProof}

\begin{Rem} \label{rem:sWc}
The time parametrization $\xi_h$ 
in the proof of Theorem~\ref{th:sWc} is nearly optimal. 
Indeed, a minimizer $\xi$ of the right hand side of \eqref{eq:WC-var} 
under the specified choice of $a$, $b$ and $\eta (r) = r$, 
as in the proof of Theorem~\ref{th:sWc}, 
is a solution to an ordinary differential equation. 
It can be expressed in the following simple form: 
\begin{align} \nonumber
\xi (r) & = l^{-1} ( \te (r) ), 
\\ \nonumber 
\te'' (r) & = - \frac{Kw^2}{2N} \mathfrak{s}_{K/N} ( 2 \te (r) ), 
\\ \nonumber
\xi (0) & = s, \, \xi (1) = t .
\end{align}
However, the solution $\te$ becomes an elliptic function in general. 
To avoid technical difficulties, 
we have considered comparison functions instead 
by linearizing the equation for $\te$ 
since we are only interested in the case $s,t \ll 1$ 
in our argument and then $\te$ must be small also. 
As a result, we obtain $\te_h$ 
in the proof of Theorem~\ref{th:sWc}. 
\end{Rem}

\section{A coupling method on Riemannian manifolds} \label{sec:coupling}

In this section, we are supposed to be in the framework of 
Theorem~\ref{th:BL-Wp}. 
In particular, $M$ is an $m$-dimensional complete Riemannian manifold 
and $P_t$ is given by the integral operator associated with 
the distribution of the diffusion process generated by 
$\sL = \Delta + Z$. 
Note that $X(t)$ is conservative under \eqref{eq:Z-CD} 
(see \cite{K10,Qian:con} for instance) 
and hence $P_t$ defines a Markov kernel. 
Since we are on a smooth space, \textbf{(A1)} is satisfied. 
% In the sequel, we study a space-time Wasserstein control 
% for $P_t$ by a coupling method based on stochastic analytic techniques. 
% Indeed, a coupling by parallel transport of two diffusion particles 
% is used to show dimension-free $W_p$-control for any $1 \le p \le \infty$ 
% (see \cite{Arn-Coul-Thal_horiz,K8,K10} and references therein). 
% To show space-time $W_p$-control, we will use a coupling by parallel transport again 
% but two particles under consideration moves in different speeds 
% to estimate distributions at different times. 

In the following argument, 
\emph{we always assume $\diam (M) < \sqrt{(N-1) \pi / K}$ when $K > 0$} 
to avoid the singularity when $d(x,y) = \sqrt{(N-1) \pi / K}$ 
(see Remark~\ref{rem:BM-BE} below). 

\begin{Rem} \label{rem:BM-BE}
When $K > 0$, it is known that 
\eqref{eq:Z-CD} yields $\diam (M) \le \sqrt{(N-1) \pi / K}$ 
(see \cite{K11}). Thus the above assumption 
only exclude the case $\diam (M) = \sqrt{(N-1)\pi/K}$. 
Even when it is the case, we can prove the same conclusion 
for $K' < K$ in \eqref{eq:Z-CD} instead of $K$ 
and finally let $K' \uparrow K$ to obtain the full statement 
from the one involving $K'$. 
As a matter of fact, $\diam (M) = \sqrt{(N-1)\pi/K}$ happens 
only when $N= \dim M$, $Z \equiv 0$ and $M$ is isometric 
to the sphere of the constant sectional curvature $K/(N-1)$ 
(see \cite{K11}). 
\end{Rem}

\subsection{The case under the absence of the cut locus} 
\label{sec:absence}

In what follows, we explain how our coupling method works. 
For this purpose, we assume that 
the cut locus of $M$ is empty and $Z = 0$ in this section 
to avoid technical difficulties. 
% On the other hand, there are some assertions 
% which we also use in the proof of the general case. 
% Indeed, it should be remarked that 
% Lemma~\ref{lem:ind_bound} and Proposition~\ref{eq:BEp} below 
% hold \emph{without this assumption}. 
In this case, we can construct a coupling of Brownian motions on $M$ 
directly by solving a coupled SDE. 
We refer to \cite{Hsu} for basic notions in this section. 

Let $\mathcal{O} (M)$ be the orthonormal frame bundle of $M$ 
and $\pi$ a canonical projection $\mathcal{O} (M) \to M$.  
Fix $\tau_2 , \tau_1 > 0$ and $x,y \in M$ with $x \neq y$ for a while. 
Let us consider a coupling of (time-scaled) horizontal diffusion processes 
$( U_1 (t) , U_2 (t) )$ by parallel transport. 
To define them, we first prepare some notations. 
For $i = 1 , \ldots , n$, 
let $H_i$ be a canonical horizontal vector field on $\mathcal{O} (M)$. 
That is, $H_i (u)$ is the horizontal lift 
(associated with the Levi-Civita connection) of 
$u e_i$, 
where $( e_i )_{i=1}^n$ is the canonical basis of $\R^n$. 
Let $\tilde{H}_i : \mathcal{O} (M) \times \mathcal{O} (M) \to T \mathcal{O} (M)$
be a horizontal vector field coupled with $( H_i )_{i=1}^n$ as follows: 
For $u_1 , u_2 \in \mathcal{O} (M)$, 
$\tilde{H}_i ( u_1 , u_2 )$ is the horizontal lift of 
$\ds 
/\!\!/_{\pi (u_1) \pi (u_2)} u_1 e_i
$, 
where $/ \!\!/_{zw}$ is the parallel transport of tangent vectors from $T_z M$ to $T_w M$
along a minimal geodesic joining $z$ and $w$ 
(such a geodesic is unique under the absence of the cut locus). 
Let $\mathbf{W} (t) = ( W^i (t) )_{i=1}^n$ be a Brownian motion on $\R^n$. 
Take $u_1 , u_2 \in \mathcal{O} (M)$ so that 
$\pi ( u_1 ) = x$, $\pi ( u_2 ) = y$ and $u_2 = /\!\!/_{xy} u_1$. 
Now we are ready to define $( U_1 (t) , U_2 (t) )$. 
They are defined as a solution to the following system of 
stochastic differential equations:  

\begin{align*}
d U_1 (t) 
& = 
\sqrt{2 \tau_1} \sum_{i=1}^n 
H_i ( U_1 (t) ) \circ d W^i (t),  
& 
U_1 (0) 
& = 
u_1 ,
\\
d U_2 (t) 
& = 
\sqrt{2\tau_2} \sum_{i=1}^n 
\tilde{H}_i ( U_1 (t), U_2 (t) ) \circ d W^i (t), 
& 
U_2 (0) 
& = 
u_2 . 
\end{align*}
Let $X_i (t) : = \pi ( U_i (t) )$. 
Then $\mathbf{X} (t) = (X_1 (t), X_2 (t) )$ is 
a coupling by parallel transport of two (time-scaled) Brownian motions. 

Intuitively, infinitesimal motions 
$d X_1 (t)$ and $d X_2 (t)$, 
which can be regarded as a random element 
in $T_{X_1 (t)} M$ and $T_{X_2 (t)} M$ respectively, 
are given by scaled ``white noises'' 
$\sqrt{ 2 \tau_1 } U_1 (t) \circ d \mathbf{W} (t)$ 
and 
$\sqrt{ 2 \tau_2 } U_2 (t) \circ d \mathbf{W} (t)$ 
respectively, and 
the second noise is given by the parallel transport of
the first noise: 
$
U_2 (t) \circ d \mathbf{W} (t) 
= 
/\!\!/_{X_1 (t) X_2 (t)} 
U_1 (t) \circ d \mathbf{W} (t)
$. 
% That is, the second noise is given by the parallel transport of the first noise. 

Let us turn to the proof. 
Let $\rho (t) : = d ( \mathbf{X} (t) )$. 
By the It\^o formula, we obtain 
\begin{align} \nonumber
d & \rho (s)^p 
=
p \rho (s)^{p-1}
\sum_{i=1}^n 
\abra{ 
  \sqrt{2 \tau_1} U_1 (s) e_i 
  \oplus  
  \sqrt{2 \tau_2} U_2 (s) e_i 
} d ( \mathbf{X} (s) ) d W^i (s) 
\\ \nonumber
& + p \rho (s)^{p-1} 
\sum_{i=1}^n 
\abra{ 
  \sqrt{\tau_1} U_1 (s) e_i 
  \oplus  
  \sqrt{\tau_2} U_2 (s) e_i 
}^2 d ( \mathbf{X} (s) ) ds 
\\ \label{eq:Ito0}
& + p (p-1) \rho(s)^{p-2}
\sum_{i=1}^n 
\bbra{
  \abra{ 
    \sqrt{\tau_1} U_1 (s) e_i 
    \oplus  
    \sqrt{\tau_2} U_2 (s) e_i 
  } d ( \mathbf{X} (s) ) 
}^2 d s .
\end{align}
We take an expectation of the integral form of \eqref{eq:Ito0}. 
To be precise, we must take care on the integrability, 
but we always assume it in this section for simpler explanation. 
The expectation of the first term in the right hand side is zero 
since it is stochastic integral. 
For the third term, by the first variation formula of 
arclength, we obtain 
\begin{align} \nonumber 
\sum_{i=1}^n 
& \bbra{
  \abra{ 
    \sqrt{\tau_1} U_1 (s) e_i 
    \oplus  
    \sqrt{\tau_2} U_2 (s) e_i 
  } d ( \mathbf{X} (s) )
}^2 
\\ \nonumber
& = 
( \sqrt{ \tau_1 } - \sqrt{ \tau_2 } )^2
\sum_{i=1}^n 
\bbra{ 
  \nabla_{U_1 (s) e_i } d ( \cdot , X_2 (s) ) 
  ( X_1 (s) )
}^2 
\\ \label{eq:bound-1}
& = 
( \sqrt{ \tau_1 } - \sqrt{ \tau_2 } )^2 . 
\end{align}

For the second term, the second variation formula of arclength, 
we obtain 

\begin{align} % \nonumber
\sum_{i=1}^n 
% & 
\abra{ 
  \sqrt{\tau_1} U_1 (s) e_i 
  \oplus  
  \sqrt{\tau_2} U_2 (s) e_i 
} d ( \mathbf{X} (s) ) 
% \\ 
\label{eq:bound0}
% & 
= 
\sum_{i=1}^n 
I_{\mathbf{X} (s)} 
\abra{ 
  \tilde{U}^{(i)}_{\mathbf{X}(s)} (s), 
  \tilde{U}^{(i)}_{\mathbf{X}(s)} (s) 
} 
ds, 
\end{align}
where $I_{zw}$ is the index form 
along a constant speed minimal geodesic from $z$ to $w$ 
and $\tilde{U}^{(i)}_{zw} (s)$ is 
the Jacobi field along the same minimal geodesic 
whose the boundary values are 
$\sqrt{ \tau_1 } U_1 (s) e_i$ and $\sqrt{ \tau_2 } U_2 (s) e_i$ 
respectively. 
For an upper bound of the index form, we introduce some notations. 
Set $K^* = K / (N-1)$. 
We define $\Psi = \Psi_{\tau_1 , \tau_2} : ( 0 , \infty ) \to \R$ 
as follows: 
\begin{equation} \label{eq:bound1}
\Psi_{\tau_1,\tau_2} (r) 
:=
(N-1) \abra{ 
  \frac{\tau_1 + \tau_2}{ \mathfrak{t}_{K^*} (r)} 
  - 
  \frac{ 2 \sqrt{ \tau_1 \tau_2 } }{ \mathfrak{s}_{K^*} (r) } 
}.
\end{equation}

\begin{Lem} \label{lem:index}
Suppose $N = m$. Let $x, y \in M$ and 
$( f_i )_{i=1}^n$ an orthonormal basis 
of $T_x M$. 
Then we have 
\begin{align*} 
\sum_{i=1}^n I_{xy} ( \tilde{V}^{(i)}_{xy} , \tilde{V}^{(i)}_{xy} ) 
& \le 
\Psi_{\tau_1,\tau_2} ( d (x,y) ),
\end{align*}
where $\tilde{V}^{(i)}_{xy}$ is the Jacobi field 
along the minimal geodesic joining $x$ and $y$ 
whose boundary values are $\sqrt{\tau_1} f_i$ and 
$\sqrt{\tau_2} /\!\!/_{xy} f_i$ respectively. 
\end{Lem}

\begin{Proof}
In the proof, we denote $d(x,y)$ by $d$. 
Recall that the condition \eqref{eq:Z-CD} is reduced 
to $\Ric \ge K$ in the present case. 
Let us define $\ph_d : [ 0 , d ] \to \R$ by 
\begin{equation} \label{eq:weight}
\ph_d (u) 
:= 
\sqrt{\tau_2} \frac{\mathfrak{s}_{K^*} (u)}{\mathfrak{s}_{K^*} ( d )} 
+ 
\sqrt{\tau_1} \frac{\mathfrak{s}_{K^*} ( d - u )}{\mathfrak{s}_{K^*} ( d )}. 
\end{equation}
Let us denote the vector field along $\gm$ 
given by the parallel transport of $f_i$ by $V_i (\cdot)$. 
Then $\ph_d (0) V_i (0) = \sqrt{\tau_1} f_i$ and 
$\ph_d (d) V_i (d) = \sqrt{\tau_2} /\!\!/_{xy} f_i$. 
Thus the index lemma together with \eqref{eq:Z-CD}
yields 
\begin{align} \nonumber
\sum_{i=1}^n I_{xy} ( V^{(i)}_{xy} , V^{(i)}_{xy} ) 
& \le 
(m-1) \int_0^{d} \ph_d' (u)^2 du 
- K \int_0^{d} \ph_d (u)^2 du 
\\ \nonumber
& = 
(m-1) \int_0^{d} \abra{ \ph_d' (u) \ph_d (u) }' du 
\\ \nonumber
& = 
(m-1) \abra{ 
  \ph_d' (d) \ph_d (d) - \ph_d' (0) \ph_d (0) 
}
\\ \nonumber
& = 
(m-1) \bbra{ 
  \frac{\sqrt{\tau_2} \mathfrak{c}_{K^*} (d) - \sqrt{\tau_1} }{ \mathfrak{s}_{K^*} (d) } \sqrt{\tau_2} 
  - 
  \frac{\sqrt{\tau_2} - \sqrt{\tau_1} \mathfrak{c}_{K^*} (d)}{\mathfrak{s}_{K^*} (d) } \sqrt{\tau_1} 
}
\\ \nonumber 
& = 
(m-1) \abra{ 
  \frac{\tau_1 + \tau_2}{ \mathfrak{t}_{K^*} (d)} 
  - 
  \frac{ 2 \sqrt{\tau_1 \tau_2} }{ \mathfrak{s}_{K^*} (d) } 
}. 
\end{align}
This is nothing but the claim. 
\end{Proof}

\begin{Lem} \label{lem:ind_bound}
For $\tau_1, \tau_2 > 0$, 
\begin{equation} \nonumber %\label{eq:bound2}
\Psi_{\tau_1,\tau_2} (r) 
\le  
\begin{cases}
\ds 
- \sqrt{ \tau_1 \tau_2 } K r
+
\frac{(N-1) ( \sqrt{ \tau_2 } - \sqrt{ \tau_1 } )^2}{r}
& (K \ge 0), 
\vspace{1.5ex}
\\ 
\ds 
- \frac{\tau_1 + \tau_2}{2} K r
+
\frac{(N-1) ( \sqrt{\tau_2} - \sqrt{\tau_1} )^2}{r} 
& ( K < 0). 
\end{cases}
\end{equation}
\end{Lem}

\begin{Proof}
By an easy rearrangement, we have 
\begin{align} \label{eq:pf-ind1}
\Psi_{\tau_1,\tau_2} (r) 
& =  
(N-1) 
\abra{ 
    2 \sqrt{\tau_1 \tau_2} \frac{\mathfrak{c}_{K^*} (r) - 1}{\mathfrak{s}_{K^*} (r)} 
    + 
    ( \sqrt{\tau_2} - \sqrt{\tau_1} )^2 \frac{1}{\mathfrak{t}_{K^*} (r)}
}. 
\end{align}
When $K \ge 0$, 
\begin{align*}
\frac{\mathfrak{c}_{K^*} (r) - 1}{\mathfrak{s}_{K^*} (r)} 
& = 
- \sqrt{K^*} \tan \abra{ \frac{\sqrt{ K^*} r}{2} } 
\le 
- \frac{K^* r}{2}, 
& 
\frac{1}{\mathfrak{t}_{K^*} (r)} 
& \le 
\frac{1}{r}.
\end{align*}
By plugging them in \eqref{eq:pf-ind1}, we obtain 
the first inequality. 
Similarly, by a rearrangement, 
\begin{align} \label{eq:pf-ind2}
\Psi_{\tau_1,\tau_2} (r)
& =  
(N-1) 
\abra{ 
    (\tau_1 + \tau_2) \frac{\mathfrak{c}_{K^*} (r) - 1}{\mathfrak{s}_{K^*} (r)} 
    + 
    ( \sqrt{\tau_2} - \sqrt{\tau_1} )^2 \frac{1}{\mathfrak{s}_{K^*} (r)}
}. 
\end{align}
When $K < 0$, 
\begin{align*}
\frac{\mathfrak{c}_{K^*} (r) - 1}{\mathfrak{s}_{K^*} (r)} 
& = 
\sqrt{-K^*} \tanh \abra{ \frac{\sqrt{-K^*}r}{2} } 
\le 
- \frac{K^* r}{2}, 
& 
\frac{1}{\mathfrak{s}_{K^*} (r)}  
& \le 
\frac{1}{r}. 
\end{align*}
and the second inequality follows 
by plugging them in \eqref{eq:pf-ind2}. 
\end{Proof}
Let us define $\tau^*$ by 
\begin{align*} % \label{eq:tau*}
\tau^*
:= 
\begin{cases}
\sqrt{\tau_1 \tau_2} 
& ( K \ge 0 ), 
\\
\ds 
\frac{\tau_1 + \tau_2}{2}
& (K < 0). 
\end{cases}
\end{align*}
By using \eqref{eq:bound-1}, \eqref{eq:bound0}, 
Lemma~\ref{lem:index} and Lemma~\ref{lem:ind_bound} 
to give an estimate of the expectation of 
the integral of \eqref{eq:Ito0} in time from $u > 0$ to $u' > u$, 
when $N = m$, we obtain 
\begin{align*}
\E [ \rho (u')^p ] 
& \le 
\E [ \rho (u)^p ] 
- p \tau^* K \int_u^{u'} \E [ \rho (v)^p ] dv
\\
& \hspace{4em} 
+ p ( N + p -2 ) 
\int_u^{u'} \E [ \rho (v)^{p-2} ] dv ( \sqrt{ \tau_1 } - \sqrt{ \tau_2 } )^2 
\\
& \le 
\E [ \rho (u)^p ] 
- p \tau^* K \int_u^{u'} \E [ \rho (v)^p ] dv
\\
& \hspace{4em} 
+ p ( N + p -2 ) 
\int_u^{u'} \E [ \rho (v)^p ]^{(p-2)/p} dv 
( \sqrt{ \tau_1 } - \sqrt{ \tau_2 } )^2 . 
\end{align*}
Obviously, the same is also true even when $N > m$.
It yields 
\begin{equation*}
\frac{\partial}{\partial u} \abra{ \E [ \rho (u)^p ]^{2/p} }
\le 
- 2 \tau^* K \E [ \rho (u)^p ]^{2/p} 
+ 2 ( N + p -2 ) ( \sqrt{ \tau_1 } - \sqrt{ \tau_2 } )^2 .
\end{equation*}
Thus, by the Gronwall lemma for $\E [ \rho (t)^p ]^{2/p}$ as a function of $t$, 
we obtain 
\begin{align*}
\E [ \rho (1)^p ]^{2/p} 
& \le 
\e^{ - 2 K \tau^* } \E [ \rho (0)^p ]^{2/p} 
+ 
\frac{(N + p - 2)( 1 - \e^{-2K \tau^* } )}{K\tau^* } 
( \sqrt{\tau_2} - \sqrt{\tau_1} )^2 .
\end{align*}
By the choice of the initial condition, 
$\E [ \rho (0)^2 ] = d ( x , y )^2$.
Moreover, since the law of $( X_1 (1) , X_2 (1) )$ 
is a coupling of 
$P_{\tau_1}^* \dl_{x}$ and $P_{\tau_2}^* \dl_{y}$, 
the definition of $L^p$-Wasserstein distance implies that 
\begin{align} \label{eq:pre-ctl}
W_p ( P_{\tau_1}^* \dl_{x} , P_{\tau_2}^* \dl_{y} )^2 
\le 
\e^{ - 2 K \tau^* } d ( x , y )^2
+ 
\frac{( N + p - 2 )( 1 - \e^{-2K \tau^* } )}{K\tau^* }  
( \sqrt{\tau_2} - \sqrt{\tau_1} )^2 .
\end{align}

To obtain \eqref{eq:Wp} from \eqref{eq:pre-ctl}, 
we prepare an $L^p$-version of the Bakry-\'Emery's curvature-dimension condition
as well as its connection with \eqref{eq:BLp}. 
For $f \in C^3 (M)$, let us define $\Gamma_2 (f)$ by 
\[
\Gamma_2 (f) = \frac12 \sL | \nabla f |^2 - \dbra{ \nabla f , \nabla \sL f }.
\]

\begin{Prop} \label{prop:BEp}
The following conditions are equivalent: 
\begin{enumerate}
\item
For any $f \in C^{\mathrm{Lip}}_b (M)$ and $x \in M$, 
\begin{align*}
| \nabla P_t f | (x)^2 
\le 
\e^{-2Kt} P_t ( | \nabla f |^{p_*} )(x)^{2/p_*}
- \frac{ 1 - \e^{-2Kt} }{ ( N + p - 2 ) K} ( \sL P_t f (x))^2 .
\end{align*}

\item
For each $t > 0$, there is a constant $C (t) > 0$ 
satisfying $\ds \lim_{t \downarrow 0} \frac{C(t)}{t} = 1$ and 
\begin{align*}
| \nabla P_t f | (x)^2 
\le 
\e^{-2Kt} P_t ( | \nabla f |^{p_*} )(x)^{2/p_*}
- \frac{2 C(t)}{N + p - 2} ( \sL P_t f (x))^2
\end{align*}
for any $f \in C^{\mathrm{Lip}}_b (M)$ and $x \in M$. 

\item
For any $f \in C^\infty (M)$, $x \in M$ and $\dl > 0$, 
\begin{align*}
\abra{ |\nabla f |^2 + \dl } 
\abra{ 
  \Gamma_2 (f) (x) 
  - K | \nabla f |(x)^2 
  - \frac{( \sL f (x))^2}{N+p-2} 
}
\ge \frac{p-2}{4(p-1)} \abs{ \nabla | \nabla f |^2 } (x)^2. 
\end{align*}
\end{enumerate}
\end{Prop}

\begin{Proof}
The implication ``(i) $\Rightarrow$ (ii)'' is obvious. 
For the proof of ``(ii) $\Rightarrow$ (iii)'', 
we claim 
\begin{equation} \label{eq:mono_app}
P_t ( ( g + \dl )^r )^{1/r} - \dl \ge P_t(  g^r )^{1/r}
\end{equation}
for $r \in ( 0 , 1 )$, $\dl > 0$ and 
$g : M \to [ 0, \infty )$ measurable. 
For \eqref{eq:mono_app}, it suffices to show 
$\partial_\dl P_t ( ( g + \dl )^r )^{1/r} \ge 1$. 
Let us take $\a > 0$ so that $\a^{-1} = r^{-1} + ( 1 - r )^{-1}$. 
Then we have 
\begin{align*}
\partial_\dl P_t ( ( g + \dl )^r )^{1/r} 
& = 
P_t ( ( g + \dl )^r )^{1/r - 1} P_t ( ( g + \dl )^{r-1} )
\\
& = 
\bbra{ 
  P_t ( ( g + \dl )^r )^{\a/r} 
  P_t \abra{ \abra{ \frac{1}{ g + \dl }}^{1-r} }^{\a/(1-r)}
}^{(1-r)/\a} 
\\
& \ge 
\bbra{ 
P_t \abra{ ( g + \dl )^\a \cdot \abra{ \frac{1}{ g + \dl } }^\a }
}^{(1-r)/\a} = 1
\end{align*}
by the H\"older inequality for $P_t$. Thus the claim holds. 

Suppose that (ii) holds. 
By applying \eqref{eq:mono_app} 
with $g = | \nabla f |^2$ and $r = p_* / 2$ 
to (ii), 
\begin{equation}
| \nabla P_t f |^2 
\le 
\e^{-2Kt} \abra{ P_t ( ( | \nabla f |^2 + \dl )^{p_*/2} )^{2/p_*} - \dl }
- \frac{2 C(t)}{N + p - 2} ( \sL P_t f )^2 . 
\end{equation}
Since the equality holds at $t = 0$ in the last inequality, 
by taking a derivative with respect to $t$ at $t = 0$, 
we obtain 
\begin{align*}
2 \dbra{ \nabla f , \nabla \sL f }
\le 
\sL | \nabla f |^2 
+ \frac{p_* - 2}{2} \cdot 
  \frac{ \abs{ \nabla | \nabla f |^2 }^2 }{ | \nabla f |^2 + \dl }
- \frac{2}{N+p-2}( \sL f )^2 - 2K | \nabla f |^2 . 
\end{align*}
By an easy rearrangement, this inequality yields (iii). 

We turn to show ``(iii) $\Rightarrow$ (i)''. 
To deal with technical difficulties, 
we consider several bounds on derivatives of $P_s f$. 
By (iii), we have $\Gamma_2 (f) \ge K | \nabla f |^2$. 
It is well-known as Bakry-\'Emery theory that 
it yields $\Ric - ( \nabla Z )^{\flat} \ge K$. 
Then we can obtain \eqref{eq:W0} and hence \eqref{eq:BE0} 
for $f \in C_b^{\mathrm{Lip}} (M)$ 
(see \cite{K9,K10} for instance). 
Note that we avoid a standard argument in Bakry-\'Emery theory 
and take a detour to obtain \eqref{eq:BE0} here in order to take the fact that 
$Z$ can be of non-gradient type into account. 
As a result of \eqref{eq:BE0}, 
$\ds \sup_{x\in M, t \le T} | \nabla P_t f | (x) < \infty$ 
for $f \in C_b^{\mathrm{Lip}} (M)$. 

Let $\dl > 0$. First we consider the case $f \in C_0^\infty (M)$. 
Let us define $F : [ 0 , t ] \times M \to \R$ by 
\[
F ( s , x ) = ( | \nabla P_{t-s} f | (x)^2 + \dl )^{p_*/2}. 
\]
Note that $F$ is bounded. Recall $p_*\ge 2$. 
To avoid technicalities on integrability, 
we will employ a stochastic analytic argument. 
Let $s_0 \in ( 0 , t )$. 
By the It\^o formula, for $s \in [ s_0 , t )$, 
\[
M^F (s) : = F ( s , X (s) ) - F ( s_0 , X (s_0) ) 
- \int_{s_0}^s 
\abra{ 
  \frac{\partial F}{\partial s} ( u , X (u) ) 
  + \sL F ( u , X (u) )
} 
du
\]
is a local martingale. 
By using (iii), we have 
\begin{align*}
\frac{\partial F}{\partial s} & (s , X (s) ) 
+ \sL F ( s , X (s) )
\\
&  = 
p_* ( | \nabla P_{t-s} f | ( X (s) )^2 + \dl )^{p_*/2 - 1} 
\abra{ 
  \Gm_2 (P_{t-s} f) 
  + \frac{p_* - 2}{4} \cdot
  \frac{\abs{ \nabla | \nabla P_{t-s} f |^2 }^2}{ | \nabla P_{t-s} f |^2 + \dl } 
}  (X(s))
\\
& \ge 
p_* ( | \nabla P_{t-s} f | ( X (s) )^2 + \dl )^{p_*/2 - 1} 
\abra{ 
  K | \nabla P_{t-s} f |^2 
  + 
  \frac{ ( P_{t-s} ( \sL f ) )^2 }{ N + p - 2 }
} (X (s)). 
\end{align*}
Note that the last term is bounded from below even when $K < 0$. 
Thus, by localizing $M^F (s)$, taking expectation and applying the Fatou lemma, 
we obtain 
\begin{align*} % \label{eq:BEp1}
\E & [ F ( s , X(s) ) ] - \E [ F ( s_0 , X (s_0) ) ] 
\\
& \ge 
p_* \int_{s_0}^s 
\E \cbra{ 
( | \nabla P_{t-u} f | ( X (u) )^2 + \dl )^{p_*/2 - 1} 
\abra{ 
  K | \nabla P_{t-u} f |^2 
  + 
  \frac{ ( P_{t-u} ( \sL f ) )^2 }{ N + p - 2 }
} ( X(u) ) 
} ds 
\\
& \ge 
p_* \int_{s_0}^s \bigg\{ K \E [ F ( u , X (u) ) ]  
+ 
\frac{( \sL P_t f )(x)^2}{
  (N+p-2) \E \cbra{ \abra{ | \nabla P_{t-u} f |(X(u))^2 + \dl }^{1 - p_* / 2} } 
}
\\
& \qquad 
- K \dl \E [ ( | \nabla P_{t-u} f | (X(u))^2 + \dl )^{p_*/2 - 1} ] 
\bigg\} du 
\\
& \ge 
p_* \int_{s_0}^s \bigg\{ K \E [ F ( u , X (u) ) ] 
+ 
\frac{( \sL P_t f )(x)^2}{N+p-2} 
\E \cbra{ \abra{ | \nabla P_{t-u} f |(X(u))^2 + \dl }^{p_* / 2} }^{1 - 2/p_*} 
\\
& \hspace{31em} 
- K \dl^{p_*/2} 
\bigg\} du .
\end{align*}
Here we have used the Schwarz inequality and 
$P_{t-s} \sL = \sL P_{t-s}$ 
in the second inequality 
and the H\"older inequality in the third inequality. 
Note that $\E [ F ( u , X(u) ) ]$ is strictly positive 
and (uniformly) continuous in $u \in [ 0, t ]$. 
Thus, 
by virtue of the mean value theorem for $r \mapsto r^{2/p_*}$ 
and the last inequality, 
for $\ep > 0$ there is a constant and $\eta > 0$ 
being independent of $s_0$ and $s$ such that 
\begin{multline*}
\e^{-2Ks} \E [ F ( s , X(s) ) ]^{2/p_*} - \e^{-2Ks_0} \E [ F ( s_0 , X(s_0) ) ]^{2/p_*} 
\\
\ge 
\abra{ 
  \frac{ ( \sL P_t f )(x)^2 }{N+p-2} - K\dl^{p_*/2} 
} 
\int_{s_0}^s \e^{-2Kr} dr 
- \ep ( s- s_0 ) 
\end{multline*}
whenever $| s - s_0 | < \eta$. 
By taking a sum and a limit, we extend the last estimate 
for $s, s_0 \in [ 0 , t ]$ with $s > s_0$. 
Then the desired estimate holds 
by substituting $s_0 = 0$ and $s = t$, 
$\ep \downarrow 0$ and $\dl \downarrow 0$.  

Finally we will show (i) for $f \in C_b^{\mathrm{Lip}} (M)$. 
The first step is to rewrite (i) 
in integral form as in the condition \textbf{(3)}. 
For $x,y \in M$ with $x \neq y$, 
let $\gm : [ 0 , 1 ] \to M$ be a minimal geodesic from $x$ to $y$.
Then, for $f \in C_0^\infty (M)$, 
$t \ge s > 0$ and $\xi (r) = rt + (1-r) s$, we have  
\begin{align*}
| P_t f (y) & - P_s f (x)|
 = 
\abs{ \int_0^1 
  \frac{ \partial }{\partial r} 
  \abra{
    P_{\xi (r)} f ( \gm (r) )
  }
dr 
}
\\
& \le 
\int_0^1 
\abra{ 
  | \nabla P_{\xi(r)} f | (\gm (r)) d (x,y)
  + 
  (t-s) \abs{ \sL P_{\xi(r)} f } ( \gm(r) )
} dr 
\\
& \le 
\int_0^1
\sqrt{ 
  \e^{-2K \xi(r)} d (x,y)^2 + \frac{(t-s)^2(N+p-2)K}{ \e^{2K \xi(r)} - 1 }
} P_{\xi(r)} ( | \nabla f |^{p_*} )( \gm (r) )^{2/p_*} dr . 
\end{align*}
Now we are ready for approximation as a second step. 
Let $f$ be Lipschitz with a compact support. 
Then there is $f_n \in C_0^\infty (M)$ 
with $\sup_n \| \nabla f_n \|_\infty < \infty$ 
such that $f_n \to f$ pointwisely and $| \nabla f_n | \to | \nabla f |$ a.e.
Since the last inequality holds for $f_n$, the same does for $f$. 
Then a usual truncation argument yields 
the same for $f \in C_b^{\mathrm{Lip}} (M)$. 
Then the conclusion follows by Proposition~\ref{prop:diff5/3}. 
\end{Proof}

\begin{Lem} \label{lem:improve}
\eqref{eq:pre-ctl} for each $x,y \in M$ and $\tau_1 , \tau_2 > 0$
implies \eqref{eq:Wp} and \eqref{eq:BLp}. 
\end{Lem} 

\begin{Proof}
By virtue of Theorem~\ref{th:Duality}, it suffices to show \eqref{eq:BLp}
for $f \in C_b^{\mathrm{Lip}} (M)$. 
By Proposition~\ref{prop:diff} and Proposition~\ref{prop:diff5/3}, 
\eqref{eq:pre-ctl} yields 
\[
| \nabla P_t f |^2 
\le 
\e^{-2Kt} P_t ( | \nabla f |^{p_*} )^{2/p_*} 
- \frac{4 K t^2 }{(N+p-2)( 1 - \e^{- 2 K t} )} ( \sL P_t f )^2 
\]
for $f \in C^{\mathrm{Lip}}_b (M)$. 
Then, by applying Proposition~\ref{prop:BEp} (ii) $\Rightarrow$ (i), 
we obtain \eqref{eq:BLp} for $f \in C^{\mathrm{Lip}}_b (M)$. 
% By an approximation, we can show the same for $f \in C_b^{\mathrm{Lip} (M)$. 
\end{Proof}

\begin{Rem} \label{rem:BEp}
It is tempting to study Proposition~\ref{prop:BEp} 
in the framework of Bakry-\'Emery theory. 
For instance, if Proposition~\ref{prop:BEp} (iii) 
follows directly from \eqref{eq:Z-CD}, then 
we can avoid the use of coupling method on which we are relying. 
However, the author do not know whether there is 
such a simpler way or not. 
Let us observe that a weaker estimate follows
by an easy application of the Bochner-Weitzenb\"ock formula:  
\[
\Gm_2 (f) 
= 
\| \Hess f \|_{\mathrm{HS}}^2
+ 
( \Ric - ( \nabla Z )^\flat )( \nabla f , \nabla f ), 
\]
where $\| \cdot \|_{\mathrm{HS}}$ stands for the Hilbert-Schmidt norm. 
Since we have 
\begin{align*}
\frac{1}{N} ( \sL f )^2 
& \le
\frac{1}{m} ( \Dl f )^2 + \frac{1}{N - m} ( Z f )^2 
\le 
\| \Hess f \|_{\mathrm{HS}}^2 + \frac{1}{N - m} ( Z f )^2 , 
\\
\abs{ \nabla | \nabla f |^2 }^2 
& \le 
4 \| \Hess f \|_{\mathrm{HS}}^2 | \nabla f |^2 , 
\end{align*}
one can show that \eqref{eq:Z-CD} yields 
\begin{align*}
( | \nabla f |^2 + \dl )
\abra{ 
  \Gm_2 (f) - K | \nabla f |^2  
  - \frac{1}{N(p-1)} ( \sL f )^2
}
& \ge 
\frac{p-2}{4(p-1)} \abs{ \nabla | \nabla f |^2 }^2, 
\\
( | \nabla f |^2 + \dl )
\abra{ 
  \Gm_2 (f) - K | \nabla f |^2  
  - \frac{1}{N+p-2} ( \sL f )^2
}
& \ge 
\frac{p-2}{4(N + p - 2)} \abs{ \nabla | \nabla f |^2 }^2 , 
\end{align*}
both of which are weaker than Proposition~\ref{prop:BEp} (iii) 
(Recall $N \ge m \ge 1$ and $p \ge 2$). 
\end{Rem}

We close this section with noting that 
Lemma~\ref{lem:index}, Lemma~\ref{lem:ind_bound}, 
% and Lemma~\ref{lem:improve}
Proposition~\ref{prop:BEp} and Lemma~\ref{lem:improve}
are all valid without the absence of the cut locus. 
Indeed, we will use them again in the next section 
(For Lemma~\ref{lem:index}, we will use a generalization of it). 

\subsection{Coupling method via discrete approximation}
\label{sec:RW}

To make the argument in the last section rigorous 
even in the presence of the cut locus, 
we will approximate the coupling of diffusion processes 
by a sequence of couplings of geodesic random walks. 
% we will construct a coupling by parallel transport of 
% time-scaled diffusion processes via approximation of 
% the diffusion process by geodesic random walks. 
% as in \cite{K8,K10,K12,Renes_poly}. 

Let $( \gm_{xy} )_{x,y \in M}$ be 
a family of unit-speed minimal geodesics 
defined on $[0, d (x,y)]$ 
such that $\gm_{xy}$ goes from $x$ to $y$. 
By using a measurable selection theorem (e.g. \cite[Theorem~6.9.6]{Boga}), 
we will take $\gm_{xy}$ as a measurable function of $(x,y)$
(more precisely, we will take a measurable choice of 
constant speed geodesics parametrized on $[ 0, 1 ]$ and 
we take $( \gm_{xy} )_{x,y}$ as their reparametrization). 
Without loss of generality, 
we may assume that 
$\gm_{xy}$ is symmetric, that is, 
$
\gm_{xy} ( d(x,y) - s ) 
= 
\gm_{yx} (s)
$ 
holds. 
Similarly as in the last section, 
we denote the parallel transport 
along $\gm_{xy}$ by $/\!\!/_{xy}$. 
We use the same symbol 
for parallel transport of orthonormal frames. 
Set $D (M) : = \{ (x,x) \; | \; x \in M \}$.  
Let $\Phi \: : \: M \to \mathscr{O} (M)$ 
be a measurable section of $\mathscr{O} (M)$. 
Let us define two measurable maps 
$\Phi_i \: : \: M \times M \to \mathscr{O} (M)$ 
for $i=1,2$ by 
\begin{align*}
\Phi_1 (x,y) & 
:= \Phi (x), 
\\ 
\Phi_2 (x,y) 
& := 
\begin{cases}
/\!\!/_{xy} \Phi_1 (x,y), 
& 
(x,y) \in M \times M \setminus D(M), 
\\ 
\Phi (x), 
& 
(x,y) \in D (M).
\end{cases}
\end{align*}
Let $( \zeta_n )_{n \in \N}$ be 
independent and identically distributed random variables 
whose distributions are uniform on the unit disk on $\R^m$. 
Take $x_1 , x_2 \in M$ and $\tau_2 > \tau_1 > 0$. 
Set $t_n^{(k)} : = k^{-2} n$ for $k \in \N$ and $n \in \N_0$ 
with $n \le k^2$. 
By using $\Phi_i$, 
we define a coupled geodesic random walk 
$\mathbf{X}^k (t) = ( X_1^k (t) , X_2^k (t) )$ 
with a discretization parameter $k \in \N$ 
by $X^k_i (0) = x_i$ and, 
for $t \in [ t_n^{(k)} , t_{n+1}^{(k)} ]$, 
\begin{align} 
\nonumber
\tilde{\zeta}_{n+1}^i
& : = 
\sqrt{2(m+2)} 
\Phi_i
\abra{ 
  \mathbf{X}^k ( t_n^{(k)} ) 
} 
\zeta_{ n + 1 } ,
\\ \nonumber 
X_i^k ( t ) 
& := 
\exp_{X_i^k ( t_n^{(k)} )}
\bigg(
  \sqrt{\tau_i}
  k^2 ( t - t_n^{(k)} )  
  \Big( 
    k^{-1} \tilde{\zeta}_{n+1}^i
     + 
    k^{-2} Z 
  \Big)
\bigg)  
\end{align}
for $i = 1, 2$, 
where $\exp_x$ is 
the exponential map at $x$. 
By \cite[Theorem~3.1]{K10} (see references therein also), 
$X_i^k (t)$ converges in law in $C ( [ 0 , \infty ) \to M )$ 
to an time-scaled $\sL$-diffusion process 
with scale parameter $\tau_i$ 
starting from $x_i$ for $i = 1,2$ respectively. 
Thus $( \mathbf{X}^k )_{k \in \N}$ is tight and 
hence a subsequential limit 
$\mathbf{X}^{k_j} \to 
\mathbf{X} = ( X_1 , X_2 )$ 
in law in $C ( [ 0, \infty ) \to M \times M )$ exists. 
Here ``time-scaled by $\tau_i$'' means that 
the law of $X_i (\cdot)$ is the same as 
$\P_{x_i} \circ X ( \tau_i \, \cdot )$. 
We fix such a subsequence $( k_j )_{j \in \N}$. 
In the rest of this paper, 
we use the same symbol $\mathbf{X}^k$ 
for the subsequence $\mathbf{X}^{k_j}$ 
and the term ``$k \to \infty$'' always means 
the subsequential limit ``$j \to \infty$''. 

Set $\rho^{(k)} (n) : = d ( \mathbf{X}^k (t_n^{(k)}) )$. 
We first show a difference inequality for $\rho^{(k)} (n)$, 
which corresponds to the It\^o formula (Lemma~\ref{lem:DI} below). 
To state it, we further introduce some notations. 
Let $\tilde{\zeta}_{n+1}^\perp (0)$ be the orthogonal projection of 
$\tilde{\zeta}_{n+1}$ to the hyperplane being perpendicular 
to $\dot{\gm}_{\mathbf{X}^k (t_n^{(k)})}(0)$. 
We denote a vector field along $\gm_{\mathbf{X}^k (t_n^{(k)})}$ given 
by parallel transport of $\tilde{\zeta}_{n+1}^\perp (0)$ 
by $( \tilde{\zeta}_{n+1}^\perp (s) )_{s \in [ 0 , \rho^{(k)} (n)]}$ and 
we define $V_{n+1} (s) := \ph_{\rho^k (n)} (s) \tilde{\zeta}_{n+1}^\perp (s)$, 
where $\ph_{\rho^k (n)}$ was defined in \eqref{eq:weight}. 
Take $v \in \R^m$ with $| v | = 1$. 
Let us define $\lm_{n+1}$ and $\Lm_{n+1}$ 
by 
\begin{align*}
\lm_{n+1} 
& : = 
\begin{cases}
\dbra{ 
  \tilde{\zeta}_{n+1}^1 (0) , 
  \dot{\gm}_{\mathbf{X}^k ( t_n^{(k)} )} (0) 
}
& \mbox{if $\mathbf{X}^k (t_n^{(k)}) \notin D(M)$}, 
\\
\sqrt{2 (m+2)} 
\dbra{ \zeta_{n+1} , v } 
& \mbox{otherwise}, 
\end{cases}
\\
\Lm_{n+1} 
& : = 
\Bigg( 
\left. 
  \ph_{\rho^k (n)} (s)
  \dbra{ Z ( t_n^{(k)} ) , \dot{\gm}_{\mathbf{X}^k (t_n^{(k)})} (s) }
\right|_{s=0}^{\rho^k (n)}
+ 
\frac12 
I_{\mathbf{X}^k (t_n^{(k)})} \abra{ V_{n+1} , V_{n+1} } 
\Bigg) 1_{\{ \mathbf{X}^k ( t_n^{(k)} ) \notin D(M) \} }, 
\end{align*}
where $I_{zw}$ stands for the index form associated with $\gm_{zw}$ 
as in Section~\ref{sec:absence}. 
Fix a reference point $o \in M$. 
For $R > 0$, let us define 
$\sg_R : C ( [ 0 , \infty ) \to M \times M ) \to [ 0, \infty ]$ 
by 
$
\sg_R ( w_1 , w_2 ) 
  : = 
\inf \bbra{ 
  t \ge 0 
  \; | \;
  d ( o , w_1 (t) ) \vee d ( o , w_2 (t) ) \ge R 
}
$. 
By \cite[Proposition~3.4]{K10}, we have the following: 
\begin{equation} \label{eq:cons}
\limsup_{R \to \infty} 
\limsup_{k \to \infty} 
\P [ \sg_R ( \mathbf{X}^k )  < \infty ] = 0. 
\end{equation}
We denote the discretization of $\sg_R ( \mathbf{X}^k )$ 
by $\hat{\sg}_R$, 
that is, 
\[
\hat{\sg}_R 
: = 
\min \{ 
  n \in \N 
  \; | \;  
  t_{n-1}^{(k)} < \sg_R ( \mathbf{X}^k ) \le t_{n}^{(k)} 
\}.
\]

\begin{Lem} \label{lem:DI}
Let $g \in C^2 ( [ 0 , \infty ) )$ be non-decreasing 
with $g' (0) = 0$ and $R > 0$ sufficiently large. 
Then there exists $k_0 \in \N$ 
such that, for any $k \ge k_0$, 
\begin{align*}
g ( \rho^{(k)} (n+1) )
& \le 
g ( \rho^{(k)} (n) )
+ \frac{1}{k} ( \sqrt{\tau_2} - \sqrt{\tau_1} ) g' ( \rho^{(k)} (n) ) \lm_{n+1} 
+ \frac{1}{k^2} g' (\rho^{(k)} (n)) \Lm_{n+1}
\\ 
& \quad 
+ \frac{1}{2k^2} ( \sqrt{\tau_2} - \sqrt{\tau_1} )^2 g'' ( \rho^{(k)} (n) ) \lm_{n+1}^2 
+ \frac{1}{R k^2} 
\end{align*}
holds on $\{ n < \hat{\sg}_{R} \} \cap \{ \rho^{(k)} (n) > R^{-1} \}$. 
\end{Lem}
\begin{Proof}
When $\Cut = \emptyset$, 
the assertion is an immediate consequence of 
the Taylor expansion, the first and second variational formulae 
and the index lemma. 
To take singularity at the cut locus into account, 
we will develop a more detailed argument 
based on the idea in \cite[Lemma~4.4]{K10}. 

Let us define $H \subset M^3$ and 
$p_1 , p_2 \: : \: H \to M^2$ 
by 
\begin{align*}
H 
& : = 
\bbra{ 
  (x,y,z) \in M^3 
  \; \left| \; 
  \begin{array}{l}
    x,y,z \in \overline{B ( o , 3 R )}, 
    \\
    d (x,y) \ge R^{-1} , 
    \\
    d (x,y) = 2 d (x,z) = 2 d (y,z) 
  \end{array}
  \right.
},
\\
p_1 & (x,y,z) 
: = 
(x,z), 
\\
p_2 & (x,y,z) 
: = 
(y,z). 
\end{align*} 
If $\mathbf{q} = (x,y,z) \in H$, 
then 
$p_1 (\mathbf{q} ) , p_2 ( \mathbf{q} ) \notin \Cut$ 
since $z$ is on a midpoint of 
a minimal geodesic joining $x$ and $y$. 
Since $H$ is compact, 
$p_1 (H)$ and $p_2 (H)$ are also compact. 
Hence both $p_1 (H)$ and $p_2 (H)$ are uniformly 
away from $\Cut \cap \overline{B ( o , 3 R )}$ 
since the cut locus is closed. 
Set 
\begin{align*}
z_n 
& := 
\gm_{\mathbf{X}^k (t_n^{(k)})} 
\abra{ 
  \frac{ \rho^{(k)} (n) }{2} 
}, 
\\
z_n'
& := 
\exp_{z_n}
\abra{ 
  V_{n+1} 
  \abra{ 
    \frac{ \rho^{(k)} (n) }{2} 
  }
  + 
  \ph_{\rho^{(k)} (n)} \abra{ \frac{\rho^{(k)}(n)}{2} }
  Z 
} .
\end{align*}
By the triangle inequality, 
we have 
\begin{align*}
\rho^{(k)} (n)
& =  
d ( X_1^k ( t_{n}^{(k)}  ) , z_n )
 + 
d ( z_n , X_2^k ( t_{n}^{(k)} ) ),    
\\
\rho^{(k)} (n+1)
& \le 
d ( X_1^k ( t_{n+1}^{(k)}  ) , z_n ' )
 + 
d ( z_n' , X_2^k ( t_{n+1}^{(k)} ) ). 
\end{align*}
Let us denote the difference of the segmented distances by $\Theta_1$ and $\Theta_2$, 
that is, 
\begin{align*}
\Theta_1 
& : = 
d ( X_1^k ( t_{n+1}^{(k)}  ) , z_n ' )
 - 
d ( X_1^k ( t_{n}^{(k)}  ) , z_n ), 
\\
\Theta_2 
& : = 
d ( z_n' , X_2^k ( t_{n+1}^{(k)} ) )
  - 
d ( z_n , X_2^k ( t_{n}^{(k)} ) ). 
\end{align*}
Suppose $t_n^{(k)} < \sg_{R} (\mathbf{X}^k )$ and 
$\rho^{(k)}(n) \ge R^{-1}$. 
Since $g$ is non-decreasing, we obtain 
\begin{align} \nonumber
g & ( \rho^{(k)} (n+1) ) - g ( \rho^{(k)} (n) ) 
\\ \nonumber
& \le 
g ( 
  d ( X_1^k ( t_{n+1}^{(k)}  ) , z_n ' )
   + 
  d ( z_n' , X_2^k ( t_{n+1}^{(k)} ) )
) 
 - 
g ( 
  d ( X_1^k ( t_{n}^{(k)}  ) , z_n )
   + 
  d ( z_n , X_2^k ( t_{n}^{(k)} ) )
)
\\ \label{eq:Taylor}
& \le 
g' ( \rho^{(k)} (n) ) ( \Theta_1 + \Theta_2 ) 
+ \frac{g''( \rho^{(k)} (n) )}{2} ( \Theta_1 + \Theta_2 )^2 
+ \frac{1}{R} ( \Theta_1 + \Theta_2 )^2
\end{align}
for sufficiently large $k$, uniformly in the position of 
$\mathbf{X}^k ( t_n^{(k)} )$.  
Note that $( X_1^k ( t_{n+1}^{(k)}  ) , z_n ' )$ is 
uniformly away from the cut locus 
since $( \mathbf{X}^k ( t_n^{(k)} ),  z_n ) \in H$. 
Thus the first and second variational formulae yield that, 
by denoting the index form along the restriction of $\gm_{\mathbf{X}^k (t_n^{(k)} )}$ 
to the geodesic from $X_1^k (t_n^{(k)})$ to $z_n$ by $I_1$, 
\begin{align*}
& \abs{ 
  \Theta_1 
   - 
  \frac{1}{k} 
  \left.
      \ph_{\rho^{(k)} (n)}(s) 
  \right|_{s=0}^{\rho^{(k)} (n)/2} 
  \lm_{n+1}
} 
\le \frac{1}{Rk}, 
\\
& 
\Theta_1 \le 
\frac{1}{k} 
\left.
    \ph_{\rho^k (n)}(s) 
\right|_{s=0}^{\rho^{(k)} (n)/2} 
\lm_{n+1}
\\
& \hspace{2em}
 + 
\frac{1}{k^2} 
\Bigg( 
\left.
    \ph_{\rho^{(k)}(n)} (s) 
    \dbra{ 
      Z , 
      \dot{\gm}_{\mathbf{X}^k (t_n^{(k)})} 
    }_{\gm_{\mathbf{X}^k (t_n^{(k)})} (s)}
\right|_{s=0}^{\rho^{(k)} (n) / 2}
+ 
\frac12 
I_{1} \abra{ V_{n+1} , V_{n+1} } 
\Bigg) 1_{\{ \mathbf{X}^k ( t_n^{(k)} ) \notin D(M) \} }
\\
& \hspace{2em} 
+ \frac{1}{R k^2}
\end{align*}
for sufficiently large $k$ uniformly in the position of 
$\mathbf{X}^k ( t_n^{(k)} )$. 
In the same way, the corresponding estimate 
also holds true for $\Theta_2$. 
Then the assertion follows 
by plugging these estimates into \eqref{eq:Taylor}. 
\end{Proof}

The next lemma estimates the expectation of 
the second variation term in Lemma~\ref{lem:DI}. 
We will use the convention 
$\sum_{n=i}^{i'} b_n = 0$ for any sequence $( b_n )_n$ 
when $i' < i$. 

\begin{Lem} \label{lem:expect1}
Let $n_1 , n_2 \in \N_0$ with $n_1 < n_2 \le k^2$ and 
$R > 0$ sufficiently large. 
Let $\mathcal{F}_n : = \sg ( \zeta_1 , \ldots , \zeta_n )$. 
\begin{enumerate}
\item \label{i:mart}
For each $h \in C ( [ 0, \infty ) \to [ 0 , \infty ) )$ and 
two $\mathcal{F}_n$-stopping times $S, T$ 
with $n_1 \wg \hat{\sg}_R \le S \le T \le n_2 \wg \hat{\sg}_R$, 
\begin{align*}
\E \cbra{ 
  \sum_{i = S + 1 }^{T}
  h ( \rho^{(k)} (i-1) ) \Lm_i
}
\le
\E \cbra{ 
  \sum_{i = S + 1 }^{T}
  h ( \rho^{(k)} (i-1) ) \Psi ( \rho^{(k)} (i-1))
}.
\end{align*}

\item \label{i:Ito1}
Let $g \in C^2 ( [ 0 , \infty ))$ be non-decreasing, $g' (0) = 0$ and 
$g'' (0) \ge 0$. 
Then there exists $k_0 \in \N$ and $C > 0$ being independent of $n_1, n_2$ and 
$R > 1$ such that, for any $k \ge k_0$, 
\begin{align*}
\E & [ g ( \rho^{(k)} ( n_2 \wg \hat{\sg}_R ) ) ] 
 \le 
\E [ g ( \rho^{(k)} ( n_1 \wg \hat{\sg}_R ) ) ]
\\
& + 
\frac{1}{k^2}
\E \Bigg[ 
  \sum_{ i = n_1 \wg \hat{\sg}_R + 1}^{ n_2 \wg \hat{\sg}_R } 
  g' ( \rho^{(k)} (i-1) ) \Psi ( \rho^{(k)} (i-1) )
   +
  g'' ( \rho^{(k)} (i-1) )   ( \sqrt{\tau_2} - \sqrt{\tau_1} )^2 
\Bigg] 
\\
& 
 + 
\frac{C ( n_2 - n_1 )}{R k^2} . 
\end{align*}
\end{enumerate}
\end{Lem}

\begin{Proof}
\ref{i:mart}
Let 
$
\bar{\Lm}_n 
 := 
\E [ 
  \Lm_n 1_{ \{ n_1 \wg \hat{\sg}_R < n \le n_2 \wg \hat{\sg}_R \} } 
   | 
  \mathcal{F}_{n-1} 
]
$. 
Then we have 
\begin{align*}
\E \cbra{ 
  \sum_{i = S + 1 }^{ T }
  h ( \rho^{(k)} (i-1) ) \Lm_i
}
= 
\E \cbra{ 
  \sum_{i = S + 1 }^{ T }
  h ( \rho^{(k)} (i-1) ) \bar{\Lm}_i
}
\end{align*}
since 
$
\sum_{i = S \wg n + 1}^{T \wg n} 
h ( \rho^{(k)} (i-1) ) ( \Lm_i - \bar{\Lm}_i )
$ 
is an $\mathcal{F}_n$-martingale. 
Let $( \hat{e}_l )_{l=1}^m$ be an orthonormal basis of $T_{X_1^k (t_n^{(k)})} M$ 
with $\hat{e}_1 = \dot{\gm}_{\mathbf{X}^k (t_n^{(k)})} (0)$ 
and $\hat{e}_l (s)$ a vector field along $\gm_{\mathbf{X}^k (t_n^{(k)})}$ 
given by parallel transport of $\hat{e}_l$. 
We also define a vector field $\bar{V}^{(l)}$ 
along $\gm_{\mathbf{X}^k (t_n^{(k)})}$ 
by $\bar{V}^{(l)}(s) := \ph_{\rho^{(k)} (n)} (s) \hat{e}_l (s)$ 
for $l = 2, \ldots , n$. 
For components of 
$\zeta_n = ( \zeta_n^{(1)} , \ldots , \zeta_n^{(m)} )$, 
we have $\E [ \zeta_n^{(l)} \zeta_n^{(l')} ] = (m+2)^{-1} \dl_{l l'}$. 
It yields 
\begin{align*}
\bar{\Lm}_n 
& = 
\Bigg( 
  \left. 
      \ph_{\rho^{(k)} (n)} (s)
      \dbra{ 
        Z , \dot{\gm}_{\mathbf{X}^k (t_n^{(k)})} 
      }_{\gm_{\mathbf{X}^k (t_n^{(k)})}(s)}
  \right|_{s=0}^{\rho^{(k)} (n)}
+ 
\sum_{l=2}^m 
I_{\mathbf{X}^k (t_n^{(k)})} \abra{ \bar{V}^{(l)} , \bar{V}^{(l)} } 
\Bigg) 
1_{\{ \mathbf{X}^k ( t_n^{(k)} ) \notin D(M) \} }
\end{align*}
on $\{ n_1 \wg \hat{\sg}_R < n \le n_2 \wg \hat{\sg}_R \}$. 
Then, based on \eqref{eq:Z-CD}, 
a similar argument as \cite[Lemma 3.4]{K12} yields 
\begin{equation*}
\bar{\Lm}_n 
\le 
(N-1) \int_0^{\rho^{(k)} (n)} \ph_{\rho^{(k)} (n)}' (u)^2 du 
- K \int_0^{\rho^{(k)} (n)} \ph_{\rho^{(k)} (n)} (u)^2 du 
\end{equation*}
on $\{ n_1 \wg \hat{\sg}_R  n \le n_2 \wg \hat{\sg}_R \}$. 
Thus the conclusion holds 
in a similar way as in Lemma~\ref{lem:index}. 

\ref{i:Ito1} 
Let us define a sequence of $\mathcal{F}_n$-stopping times $S_j$ 
($j = 0 , 1 , \ldots $) as follows: 
\begin{align*}
S_0 & : = n_1 \wg \hat{\sg}_R , 
\\
S_{2j+1} & : = 
\min \{ n \ge S_{2j} \; | \; \rho^{(k)} (n) < 2 R^{-1} \} 
\wg n_2 \wg \hat{\sg}_R , 
\\
S_{2j+2} & : = 
\min \{ n > S_{2j+1} \; | \; \rho^{(k)} (n) > \rho^{(k)} ( S_{2j+1} ) \} 
\wg n_2 \wg \hat{\sg}_R . 
\end{align*}
For simplicity of notations, set 
\[
L (i) : = 
  g' ( \rho^{(k)} (i-1) ) \Psi ( \rho^{(k)} (i-1) )
   +
  g'' ( \rho^{(k)} (i-1) )   ( \sqrt{\tau_2} - \sqrt{\tau_1} )^2 . 
\]
Note that 
$\E [ \lm_{n+1} | \mathcal{F}_n ] = 0$ and 
$\E [ \lm_{n+1}^2 | \mathcal{F}_n ] = 2$. 
Thus, as an immediate consequence of the first assertion and Lemma~\ref{lem:DI}, 
\begin{align} \label{eq:Ito1}
\E 
[ g ( \rho^{(k)} ( S_{2j+1} ) ) ] 
\le 
\E [ g ( \rho^{(k)} ( S_{2j} ) ) ]
+ 
\frac{1}{k^2}
\E \Bigg[ 
  \sum_{ i = S_{2j} + 1}^{ S_{2j+1} } 
  L(i)
\Bigg] 
 + 
\frac{ \E \cbra{ S_{2j+1} - S_{2j} }}{R k^2}  
\end{align}
when $k$ is sufficiently large. 
On the other hand, since $\rho^{(k)} ( S_{2j+2} - 1 ) > R^{-1}$ 
if $S_{2j+1} < n_2 \wg \hat{\sg}_R$ and $k$ is sufficiently large, 
\begin{align} \nonumber
\E [ g ( \rho^{(k)} 
( S_{2j+2} ) ) ] 
& \le 
\E [ g ( \rho^{(k)} ( S_{2j+1} ) ) ]
\\ \nonumber
& \qquad 
+ \E [ ( g ( \rho^{(k)} ( S_{2j+2} ) )  
-  g ( \rho^{(k)} ( S_{2j+2} - 1 ) ) ) 1_{\{ S_{2j+1} < n_2 \wg \hat{\sg}_R \}} ] 
\\ \nonumber
& \le \E [ g ( \rho^{(k)} ( S_{2j+1} ) ) ]
\\ \label{eq:Ito2}
& \qquad + 
\frac{1}{k^2}
\E \Bigg[ 
L ( S_{2j+2} ) 
1_{\{ S_{2j+1} < n_2 \wg \hat{\sg}_R \}}
\Bigg]
 + 
\frac{ 
  \P \cbra{ S_{2j+1} < n_2 \wg \hat{\sg}_R 
  }
}{R k^2} .  
\end{align}
By using $g' (0) = 0$ and \eqref{eq:pf-ind1}, 
we can show that there is a function $c : ( 0 , \infty ) \to \R$ 
with $\lim_{r \to 0} c(r) = 0$ 
such that $g' (r) \Psi (r) \ge c (r)$. 
This fact together with $g'' (0) \ge 0$ implies that 
there is a constant $C > 0$ such that 
\begin{align} \label{eq:Ito3}
0 \le 
\frac{1}{k^2}
\E \Bigg[ 
  \sum_{ i = S_{2j+1} + 1}^{ S_{2j+2} - 1 } 
  L (i)
\Bigg] 
 + 
\frac{C}{R k^2}
  \E \cbra{ 
     ( S_{2j+2} - S_{2j+1} - 1 ) \vee 0 
  }. 
\end{align}
Note that $C$ can be chosen to be independent of $n_1$ and $n_2$. 
$C$ may depend on $R$ but it can be smaller for larger $R$. 
Thus we can choose it to be independent of $R$ also. 
Then the assertion holds  
by summing up \eqref{eq:Ito1}, \eqref{eq:Ito2} and \eqref{eq:Ito3} 
and take a summation again with respect to $j$. 
\end{Proof}

The third lemma deals with the limit $k \to \infty$ and 
a Gronwall type bound for expectations for truncated functions. 

\begin{Lem} \label{lem:limit}
Let $g$ be as in Lemma~\ref{lem:expect1} \ref{i:Ito1}. 
Let $\psi \in C^2 ( [ 0,  \infty ) )$ be 
an increasing concave function satisfying $\psi (x) = x$ 
for $x \in [ 0, 1 ]$ and $\psi (x) = 2$ for $x \in [ 3 , \infty )$ 
and set $\psi_j (x) : = j \psi ( x / j )$ and $g_j := g \circ \psi_j$ 
for $j \in \N$. 
Then, for $0 \le s_1 < s_2 \le 1$, 
\begin{align*}
\E [ g_j ( d ( \mathbf{X} ( s_2 ) ) ) ] 
& \le 
\E [ g_j ( d ( \mathbf{X} ( s_1 ) ) ) ]
+ 
\int_{s_1}^{s_2} 
\E \cbra{ 
  ( g_j' \cdot \Psi ) ( d ( \mathbf{X} ( u ) ) )
}
d u 
\\ %\label{eq:pre-Gron2}
& \hspace{7em}
+
( \sqrt{\tau_2} - \sqrt{\tau_1} )^2 
\int_{s_1}^{s_2} 
\E \cbra{ g_j'' ( d ( \mathbf{X} ( u ) ) ) }
d u . 
\end{align*}
\end{Lem}

\begin{Proof}
Take $R > 0$ sufficiently large. 
We set 
$
\lfloor s_i \rfloor_k 
:= 
\inf \{ n \in \N \; | \; t_n^{(k)} < s_i \le t_{n+1}^{(k)} \} 
$.
Note that $g_j$, $g_j' \cdot \Psi$ and $g_j''$ are 
all uniformly continuous and bounded 
on $[ 0, \infty )$. 
Thus, by Lemma~\ref{lem:expect1} \ref{i:Ito1}, 
for sufficiently large $k$, we have 
\begin{align} \nonumber
\E & [ g_j ( d ( \mathbf{X}^k ( s_2 \wg \sg_R ) ) ) ] 
 \le 
\E [ g_j ( d ( \mathbf{X}^k ( s_1 \wg \sg_R ) ) ) ]
\\ \nonumber
& \quad + 
\frac{1}{k^2}
\E \Bigg[ 
  \sum_{ i = \lfloor s_1 \rfloor_k \wg \hat{\sg}_R + 1}^{ \lfloor s_2 \rfloor_k \wg \hat{\sg}_R } 
  ( g_j' \cdot \Psi ) ( d ( \mathbf{X}^{k} ( t^{(k)}_{i-1} ) ) )
   +
  g_j'' ( d ( \mathbf{X}^k ( t^{(k)}_{i-1} ) ) )
  ( \sqrt{\tau_2} - \sqrt{\tau_1} )^2 
\Bigg] 
 + 
\frac{2}{R} 
\\ \nonumber
& \le 
\E [ g_j ( d ( \mathbf{X}^k ( s_1 \wg \sg_R ) ) ) ]
+ 
\E \cbra{ 
  \int_{s_1 \wg \sg_R}^{s_2 \wg \sg_R} 
  ( g_j' \cdot \Psi ) ( d ( \mathbf{X}^{k} ( u ) ) )
  d u 
}
\\ \label{eq:pre-Gron1}
& \hspace{2em} +
( \sqrt{\tau_2} - \sqrt{\tau_1} )^2 
\E \cbra{
  \int_{s_1 \wg \sg_R}^{s_2 \wg \sg_R} 
  g_j'' ( d ( \mathbf{X}^k ( u ) ) )
  d u 
}
+ \frac{3}{R}. 
\end{align}
Since 
$\{ \mathbf{w} \; | \; \sg_{R-1} ( \mathbf{w} ) > s \}$ is open 
in $C ( [ 0 , \infty ) \to M \times M )$, 
the Portmanteau theorem 
for the weak convergence $\mathbf{X}^k \to \mathbf{X}$ 
yields 
\begin{multline*}
\liminf_{R \to \infty} 
\liminf_{k \to \infty} 
\E [ g_j ( d ( \mathbf{X}^k ( s_2 \wg \sg_R ) ) ) ] 
\ge 
\liminf_{R \to \infty} 
\liminf_{k \to \infty} 
\E [ g_j ( d ( \mathbf{X}^{k} ( s_2 ) ) ) \; ; \; \sg_{R} ( \mathbf{X}^k ) > s_2 ] 
\\
\ge 
\liminf_{R \to \infty} 
\E [ g_j ( d ( \mathbf{X} ( s_2 ) ) ) \; ; \; \sg_{R} ( \mathbf{X} ) > s_2 ] 
= 
\E [ g_j ( d ( \mathbf{X} ( s_2 ) ) ) ], 
\end{multline*}
where the last inequality follows 
from the fact that $X$ is conservative. 
On the other hand, for any $h \in C_b ( [ 0, \infty ) )$, 
\eqref{eq:cons} yields 
\begin{align*}
\limsup_{R \to \infty} 
& 
\limsup_{k \to \infty} 
\E \cbra{ 
  \int_{s_1 \wg \sg_R}^{s_2 \wg \sg_R} 
    h ( d ( \mathbf{X}^k ( u ) ) ) 
  d u 
}
\\
& \le
\| h \|_\infty ( s_2 - s_1 )
\limsup_{R\to\infty} \limsup_{k \to \infty} 
\P [ \sg_R ( \mathbf{X}^k ) \le s_2 ] 
\\
& \qquad + 
\limsup_{R \to \infty} \limsup_{k \to \infty} 
\E \cbra{ 
  \int_{s_1}^{s_2} 
    h ( d ( \mathbf{X}^k ( u ) ) )
  d u 
  \; ; \; 
  \sg_R ( \mathbf{X}^k ) > s_2  
}
\\
& \le 
\limsup_{k \to \infty} 
\E \cbra{ 
  \int_{s_1}^{s_2} 
    h ( d ( \mathbf{X}^k ( u ) ) ) 
  d u 
} 
= 
\int_{s_1}^{s_2} 
\E \cbra{  
  h ( d ( \mathbf{X} (u) ) ) 
}
d u . 
\end{align*}
A similar argument also works for the first term 
in the right hand side of \eqref{eq:pre-Gron1}. 
Then the conclusion follows 
by applying these estimates to \eqref{eq:pre-Gron1}, 
when we take the limit $k \to \infty$ and $R \to \infty$ after it.  
\end{Proof}

We are now ready to show 
the key assertion \eqref{eq:pre-ctl} 
in the last section. 

\begin{Prop} \label{prop:Lp}
\eqref{eq:pre-ctl} holds. 
\end{Prop}

\begin{Proof}
Let $\psi_j$ be as in Lemma~\ref{lem:limit}. 
Note first that 
we have $\psi'_j (u) u \le \psi_j (u)$ and $\psi_j'' (u) \le 0$ 
for each $u \ge 0$ since $\psi$ is concave and $\psi (0) = 0$. 
In addition, $\psi_j$ and $\psi_j'$ is non-decreasing in $j$ 
and we have 
$\lim_{j \to \infty} \psi_j (u) = u$, 
$\lim_{j \to \infty} \psi_j' (u) = 1$ and 
$\lim_{j \to \infty} \psi_j'' (u) = 0$ 
for each $u \ge 0$. 
By applying Lemma~\ref{lem:limit} with $g (u) := u^p$ 
together with Lemma~\ref{lem:ind_bound}, we obtain 
\begin{align} \nonumber
\E & [ \psi_j ( d ( \mathbf{X} (s_2 ) ) )^p ] 
\le 
\E [ \psi_j ( d ( \mathbf{X} ( s_1 ) ) )^p ]
\\ \nonumber
& \qquad -  
p K \tau^*
\int_{s_1}^{s_2} 
\E \cbra{ 
  \psi_j ( d ( \mathbf{X} ( u ) ) )^{p-1} 
  \psi_j' ( d ( \mathbf{X} (u) ) )
  d ( \mathbf{X} (u) ) 
}
d u 
\\ \nonumber
& \qquad +   
p (N-1) ( \sqrt{ \tau_2 } - \sqrt{ \tau_1 } )^2 
\int_{s_1}^{s_2} 
\E \cbra{ 
  \psi_j ( d ( \mathbf{X} ( u ) ) )^{p-1} 
  \frac{\psi_j' ( d ( \mathbf{X} (u) ) )}
  { d ( \mathbf{X} (u) ) }
}
d u 
\\ \nonumber 
& \qquad 
+
p ( p - 1 ) ( \sqrt{\tau_2} - \sqrt{\tau_1} )^2 
\int_{s_1}^{s_2} 
\E \cbra{ 
  \psi_j ( d ( \mathbf{X} ( u ) ) )^{p-2} 
  \psi_j' ( d ( \mathbf{X} ( u ) ) )^2 
} 
d u 
\\ \label{eq:pre-Gron2}
& \qquad + 
p ( \sqrt{\tau_2} - \sqrt{\tau_1} )^2 
\int_{s_1}^{s_2} 
\E \cbra{ 
  \psi_j ( d ( \mathbf{X} ( u ) ) )^{p-1} 
  \psi_j'' ( d ( \mathbf{X} (u) ) ) 
}
d u . 
\end{align}
By neglecting non-positive terms, 
trivial bounds $\psi_j (u) \le u$ and $\psi_j' (u) \le 1$, 
properties of $\psi_j$ stated at the beginning and  
the H\"older inequality yields 
\begin{align} \nonumber
\E [ \psi_j ( d ( \mathbf{X} (s_2 ) ) )^p ] 
& \le 
\E [ \psi_j ( d ( \mathbf{X} ( s_1 ) ) )^p ]
-  
p ( K \wg 0 ) \tau^*
\int_{s_1}^{s_2} 
\E [ \psi_j ( d ( \mathbf{X} ( u ) ) )^p ]
d u 
\\ \nonumber 
& \hspace{2em} +   
p ( N + p - 2 ) ( \sqrt{ \tau_2 } - \sqrt{ \tau_1 } )^2 
\int_{s_1}^{s_2} 
\E [ \psi_j ( d ( \mathbf{X} ( u ) ) )^{p-2} ]
d u 
\\ \nonumber 
& \le 
\E [ \psi_j ( d ( \mathbf{X} ( s_1 ) ) )^p ]
-  
p ( K \wg 0 ) \tau^*
\int_{s_1}^{s_2} 
\E [ \psi_j ( d ( \mathbf{X} ( u ) ) )^p ]
d u 
\\ \nonumber
& \hspace{2em} +   
p ( N + p - 2 ) ( \sqrt{ \tau_2 } - \sqrt{ \tau_1 } )^2 
\int_{s_1}^{s_2} 
\E [ \psi_j ( d ( \mathbf{X} ( u ) ) )^p ]^{1-2/p}
d u . 
\end{align}
Set $a_j (s) : = \E [ \psi_j ( d ( \mathbf{X} (s) ) )^p ]$. 
Then the last inequality implies that, for $\dl \in ( 0, 1 )$, 
\begin{align} \nonumber
\limsup_{h \to 0} & 
\frac{ 
  \e^{2 ( K \wg 0 ) \tau^* (s+h)} ( a_j (s+h) + \dl )^{2/p} 
  - \e^{2 ( K \wg 0 ) \tau^* s} ( a_j (s) + \dl )^{2/p} }{h}
\\ \nonumber
& \le 
2 ( K \wg 0 ) \tau^* \e^{2 ( K \wg 0 ) \tau^* s} 
( a_j (s) + \dl )^{2/p - 1} \dl 
\\ \nonumber 
& \qquad +
2 ( N + p - 2 ) ( \sqrt{ \tau_2 } - \sqrt{ \tau_1 } )^2 
\e^{2 ( K \wg 0 ) \tau^* s} 
( a_j (s) + \dl )^{2/p - 1} a_j (s)^{1-2/p} 
\\ \label{eq:pre-Gron3}
& \le 
2 ( N + p - 2 ) ( \sqrt{ \tau_2 } - \sqrt{ \tau_1 } )^2 
\e^{2 ( K \wg 0 ) \tau^* s} . 
\end{align}
Thus the Gronwall lemma yields that 
there is a constant $C_1 > 0$ being independent of $j,s$ and $\dl$ 
such that $a_j (s) \le C_1$ holds. 
% \begin{align} \label{eq:pre-Gron4}
% a_j (s) 
% \le 
% \e^{-2 ( K \wg 0 ) \tau^* s } a_j (0) 
% + 
% 2 ( N + p - 2 ) ( \sqrt{ \tau_2 } - \sqrt{ \tau_1 } )^2 
% \e^{-2 ( K \wg 0 ) \tau^* s } \int_0^s \e^{2 ( K \wg 0 ) \tau^* u} d u. 
% \end{align}
Therefore the monotone limit $a_\infty : = \lim_{j \to\infty} a_j (s)$ 
exists in $\R$. 
In addition, the monotone convergence theorem yields 
$a_\infty (s) = \E [ d ( \mathbf{X} (s) )^p ]$. 

With the aid of the monotone convergence theorem 
and the dominated convergence theorem, 
by letting $j \to \infty$ in \eqref{eq:pre-Gron2} 
and by applying the H\"older inequality, we obtain 
\begin{align*}
a_\infty (s_2) 
\le 
a_\infty (s_1) 
- p K \tau^* \int_{s_1}^{s_2} a_\infty (u) d u 
+ p ( N + p - 2 ) ( \sqrt{ \tau_2 } - \sqrt{ \tau_1 } )^2 
\int_{s_1}^{s_2} a_\infty (u)^{1-2/p} d u . 
\end{align*}
Then we argue as in \eqref{eq:pre-Gron3} with this inequality 
to apply the Gronwall inequality. 
Consequently, we obtain 
\begin{align*}
( a_\infty (1) + \dl )^{2/p}
& \le 
\e^{- 2 K \tau^*} ( a_\infty (0) + \dl )^{2/p}
+ 
2 ( N + p - 2 ) ( \sqrt{ \tau_2 } - \sqrt{ \tau_1 } )^2 
\e^{- 2 K \tau^*} 
\int_0^1 
\e^{2 K \tau^* u} 
d u 
\\
& \qquad + 2 ( K \vee 0 ) \e^{- 2 K \tau^*} 
\int_0^1 
\e^{2 K \tau^* u} 
d u \, 
\dl^{2/p} . 
\end{align*}
Since $a_\infty (0) = d (x,y)^p$ and 
the definition of the Wasserstein distance and $\mathbf{X}$ 
implies $a_\infty (1) \ge W_p ( P_{\tau_1}^* \dl_x , P_{\tau_2}^* \dl_y )^p$,
the conclusion holds by letting $\dl \downarrow 0$ 
in the last inequality. 
\end{Proof}

\begin{tProof}{Theorem~\ref{th:BL-Wp}}
The combination of Proposition~\ref{prop:Lp} and 
Lemma~\ref{lem:improve} immediately completes the proof. 
\end{tProof}

For $c : M \times M \to \R$ measurable and 
bounded from below, and $\mu_1 , \mu_2 \in \sP (M)$, 
we define the optimal transportation cost $\mathcal{T}_c ( \mu_1 , \mu_2 )$ 
between $\mu_1$ and $\mu_2$ associated with the cost function $c$ 
as follows: 
\[
\mathcal{T}_c ( \mu_1 , \mu_2 ) 
= 
\inf \bbra{ 
  \left.
      \int_{M \times M} c \, d \pi 
  \; \right| \; 
  \mbox{$\pi$ is a coupling of $\mu_1$ and $\mu_2$}
}. 
\]
As a variant of Proposition~\ref{prop:Lp}, we obtain the following: 

\begin{Thm} \label{th:Lp2}
For $\tau_1 , \tau_2 > 0$, $\mu_1 , \mu_2 \in \sP (M)$ 
and $p \ge 2$, 
\begin{align*}
\mathcal{T}_{\mathfrak{s}_{K^*}^p (d/2)}
( P^*_{\tau_1} \mu_1 , P^* _{\tau_2} \mu_2 )^{2/p}
% & 
\le 
\e^{-\te}
\mathcal{T}_{\mathfrak{s}_{K^*}^p (d/2)} ( \mu_1 , \mu_2 )^{2/p}
% \\
% & \qquad 
+ 
\frac{
  ( N + p - 2 ) 
  ( 1 - \e^{ - \te } )
}
{ 
  2 \te 
}
( \sqrt{\tau_2} - \sqrt{\tau_1} )^2 , 
\end{align*}
where 
$ \te = 
\te ( \tau_1 ,\tau_2 , K , N , p ) 
 := 
K ( \tau_1 + \tau_2 ) 
+ 
pK^* ( \sqrt{\tau_2} - \sqrt{\tau_1} )^2 / 2
$, 
and $K^* = K / ( N - 1 )$ 
as in the definition of $\Psi$ in \eqref{eq:bound1}.  
\end{Thm}

Before entering the proof, we recall the following 
elementary relations for comparison functions: 
\begin{align*} 
\mathfrak{s}_{K^*}' & = \mathfrak{c}_{K^*}, 
& 
\mathfrak{c}_{K^*}' & = - K^* \mathfrak{s}_{K^*}, 
& 
\mathfrak{c}_{K^*}^2 + K^* \mathfrak{s}_{K^*}^2 & = 1,
\\
\mathfrak{s}_{K^*} (2r) & = 2 \mathfrak{s}_{K^*} (r) \mathfrak{c}_{K^*} (r), 
& 
\mathfrak{c}_{K^*} (2r) & = \mathfrak{c}_{K^*} (r)^2 - K^* \mathfrak{s}_{K^*} (r)^2 .
\end{align*}

\begin{Proof}
The same argument as in Lemma~\ref{lem:Dirac} works 
for the transportation cost $\mathcal{T}_{\mathfrak{s}_{K^*}^p (d/2)}$ 
instead of $W_p^p$. 
Hence it suffices to show the assertion 
only when both $\mu_1$ and $\mu_2$ are Dirac measures. 
We consider only the case $K \neq 0$ 
since the assertion is reduced to \eqref{eq:pre-ctl} 
when $K=0$. 
We begin with the integrability of 
$\mathfrak{s}_{K^*} ( d ( \mathbf{X} (s) ) )^p$ 
as in the proof of Proposition~\ref{prop:Lp}. 
Since it is obvious when $K > 0$, 
we assume $K < 0$ for a while. 
By applying Lemma~\ref{lem:limit} 
with $g (u) = \mathfrak{s}_{K^*} ( u / 2 )^p$ 
and \eqref{eq:pf-ind2} yield 
\begin{align} \nonumber
\E & \cbra{ 
\mathfrak{s}_{K^*} \abra{
 \frac{\psi_j ( d ( \mathbf{X} (s_2) ) )}{2}
}^p 
} 
\le
\E \cbra{ 
  \mathfrak{s}_{K^*} \abra{ 
    \frac{\psi_j ( d ( \mathbf{X} (s_1) ) )}{2}
  }^p 
} 
\\ \nonumber
& 
- 
\frac{ p K ( \tau_1 + \tau_2 ) }{2}
\int_{s_1}^{s_2} 
\E \cbra{ 
  \abra{ 
    \mathfrak{s}_{K^*}^{p-1} \cdot \mathfrak{c}_{K^*} 
  }
  \abra{ 
    \frac{\psi_j ( d ( \mathbf{X} (u) ) )}{2}
  } 
  \psi_j' ( d ( \mathbf{X} (u) ) )
  \mathfrak{t}_{K^*} 
  \abra{ 
    \frac{ d ( \mathbf{X} (u) ) }{2}
  }
} 
d u 
\\ \nonumber
& + 
\frac{ p (N-1) ( \sqrt{\tau_2} - \sqrt{\tau_1} )^2 }{2}
\int_{s_1}^{s_2} 
\E \cbra{ 
  \abra{ 
    \mathfrak{s}_{K^*}^{p-1} \cdot \mathfrak{c}_{K^*} 
  }
  \abra{ 
    \frac{\psi_j ( d ( \mathbf{X} (u) ) )}{2}
  } 
  \frac{ \psi_j' ( d ( \mathbf{X} (u) ) ) }
  { \mathfrak{s}_{K^*} \abra{ d ( \mathbf{X} (u) ) } }
} 
d u 
\\ \nonumber
& + 
\frac{p ( \sqrt{\tau_2} - \sqrt{\tau_1} )^2 }{4}
\int_{s_1}^{s_2} 
\E \cbra{ 
  \abra{ 
    (p-1) \mathfrak{s}_{K^*}^{p-2} 
    - p K^* 
    \mathfrak{s}_{K^*}^p 
  } 
  \abra{ 
    \frac{ \psi_j ( d ( \mathbf{X} (u) ) ) }{2} 
  } 
  \psi_j' \abra{ d ( \mathbf{X} (u) ) }^2 
}
d u 
\\ \label{eq:s-Gron1}
& + 
\frac{ p ( \sqrt{\tau_2} - \sqrt{\tau_1} )^2 }{2}
\int_{s_1}^{s_2} 
\E \cbra{ 
  \abra{ 
    \mathfrak{s}_{K^*}^{p-1} \cdot \mathfrak{c}_{K^*} 
  }
  \abra{ 
    \frac{\psi_j ( d ( \mathbf{X} (u) ) )}{2}
  } 
  \psi_j'' ( d ( \mathbf{X} (u) ) ) 
} 
d u . 
\end{align}
Note that there exists $c_j \ge 1$ with $\lim_{j \to \infty} c_j = 1$ 
such that 
$\mathfrak{c}_{K^*} (\psi_j (u) / 2 ) \cdot \mathfrak{t}_{K^*} (u/2) 
\le c_j \mathfrak{s}_{K^*} ( \psi_j (u) / 2 )$ for each $u \ge 0$. 
Set $\tilde{a}_j (s) := \E [ \mathfrak{s}_{K^*} ( \psi_j ( d ( \mathbf{X} (s) ) ) / 2 )^p ]$. 
By a similar argument as in the proof of Proposition~\ref{prop:Lp}, 
from \eqref{eq:s-Gron1}, we obtain 
\begin{align*}
\tilde{a}_j (s_2) & \le
\tilde{a}_j (s_1)  
- \frac{c_j p \theta}{2} 
\int_{s_1}^{s_2} \tilde{a}_j (u) d u 
+ 
\frac{ p ( N + p - 2 ) ( \sqrt{\tau_2} - \sqrt{\tau_1} )^2 }{4}
\int_{s_1}^{s_2} 
\tilde{a}_j (u)^{1-2/p} d u . 
\end{align*}
Thus, as we discussed in \eqref{eq:pre-Gron3}, 
we can show that there exists $\tilde{C}_1 > 0$ 
being independent of $s$ and $j$ 
such that 
$\tilde{a}_j (s) \le \tilde{C}_1$. 
It ensures 
$\tilde{a}_\infty (s) := \lim_{j \to \infty} \tilde{a}_j (s) < \infty$ 
for each $s \ge 0$. 

Now we turn to the general situation $K \in \R$. 
\eqref{eq:s-Gron1} is still valid in this case. 
Then, by taking the limit $j \to \infty$, 
the conclusion follows in the same way 
as in the proof of Proposition~\ref{prop:Lp}. 
\end{Proof}

% \bibliographystyle{amsplain}
% \bibliography{refs}

\begin{thebibliography}{10}

\bibitem{Ambrosio:2012tp}
L.~Ambrosio, N.~Gigli, A.~Mondino, and T.~Rajala, \emph{{Riemannian Ricci
  curvature lower bounds in metric measure spaces with ${\sigma}$-finite
  measure}}, To appear in Trans. Amer. Math. Soc., 2012.

\bibitem{AGS_BE-CD}
L.~Ambrosio, N.~Gigli, and G.~Savar\'e, \emph{Bakry-{{\'E}}mery
  curvature-dimension condition and {R}iemannian {R}icci curvature bounds},
  Preprint. Available at: arXiv: 1209.5786.

\bibitem{AGS2}
\bysame, \emph{Calculus and heat flow in
  metric measure spaces and applications to spaces with {R}icci bounds from
  below}, To appear in Invent. Math. 
% Preprint. 
Available at: arXiv:1106.2090.

\bibitem{AGS3}
\bysame, \emph{Metric measure spaces with {R}iemannian {R}icci curvature
  bounded from below}, To appear in Duke Math. J. 
% Preprint. 
Available at: arXiv:1109.0222.

\bibitem{AGS}
\bysame, \emph{Gradient flows in metric spaces and in the space of probability
  measures}, second ed., Birkh{\"a}user Verlag, Basel, 2008.

\bibitem{AGS_S}
\bysame, \emph{Density of Lipschitz functions and equivalence of 
weak gradients in metric measure spaces} , Rev. Mat. Iberoam. \textbf{29}
(2013), no.~3, 969--996. 

\bibitem{Arn-Coul-Thal_horiz}
M.~Arnaudon, K.A. Coulibaly, and A.~Thalmaier, \emph{Horizontal diffusion in
  {$C^1$}-path space}, S\'{e}minaire de Probabilit\'{e}s, XLIII, Lecture Notes
  in Mathematics, 2006, Springer, Berlin, 2011, pp.~73--94.

\bibitem{Bacher:2010bp}
K.~Bacher and K.-T. Sturm, \emph{{Localization and tensorization properties of
  the curvature-dimension condition for metric measure spaces}}, J. Funct.
  Anal. \textbf{259} (2010), no.~1, 28--56.

\bibitem{Bak97}
D.~Bakry, \emph{On {S}obolev and logarithmic {S}obolev inequalities for
  {M}arkov semigroups}, New trends in stochastic analysis (Charingworth, 1994),
  World Sci. Publ. River Edge, NJ, 1997, pp.~43--75.

\bibitem{Bakry:2006ul}
\bysame, \emph{{Functional inequalities for Markov semigroups}}, Probability
  measures on groups: recent directions and trends (Mumbai), Tata Inst. Fund.
  Res., 2006, pp.~91--147.

\bibitem{BE:diff-hyp}
D.~Bakry and M.~{\'E}mery, \emph{Diffusions hypercontractives}, S{\'e}minaire
  de probabilit{\'e}s, XIX, 1983/1984, Lecture notes in Mathematics, 1123,
  Springer, Berlin, 1985, pp.~177--206.

\bibitem{BGL_book}
D.~Bakry, I.~Gentil, and M.~Ledoux, \emph{Analysis and geometry of Markov diffusion operators}
Grundlehren der Mathematischen Wissenschaften, vol.~348, Springer-Verlag, 2014.

\bibitem{BGL}
D.~Bakry, I.~Gentil, and M.~Ledoux, \emph{On {H}arnack inequalities and optimal
  transport}, To appear in Ann. Sc. Norm. Super. Pisa. 
% Preprint. 
Available at: arXiv:1210.4650.

\bibitem{BL}
D.~Bakry and M.~Ledoux \emph{Sobolev inequalities and Myers's diameter theorem 
for an abstract Markov generator}, Duke Math. J., \textbf{85} (2006) 253–270.

\bibitem{Bakry:2006vh}
D.~Bakry and M.~Ledoux, \emph{{A logarithmic Sobolev form of the Li-Yau
  parabolic inequality}}, Rev. Mat. Iberoam. \textbf{22} (2006), no.~2,
  683--702.

\bibitem{BEHM_HJ}
Z.M.~Balogh, A.~Engoulatov, L.~Hunziker, and O.E. Maasalo, \emph{Functional
  inequalities and {H}amilton-{J}acobi equations in geodesic spaces}, Potential
  Anal. \textbf{36} (2012), no.~2, 317--337.

\bibitem{Boga}
V.~I.~Bogachev, \emph{Measure theory. Vol. I, II}
Springer, Berlin, 2007. 

\bibitem{Bolley:2013ts}
F.~Bolley, I.~Gentil, and A.~Guillin, \emph{{Dimensional contraction via Markov
  transportation distance}}, To appear in J. Lond. Math. Soc. 
Available at: arXiv:1304.1929.
% Preprint. 

\bibitem{Crans}
M.~Cranston, \emph{Gradient estimates on manifolds using coupling}, J. Funct.
  Anal. \textbf{99} (1991), no.~1, 110--124.

\bibitem{Erbar:2013wf}
M.~Erbar, K.~Kuwada, and K.-T. Sturm, \emph{{On the equivalence of the entropic
  curvature-dimension condition and Bochner's inequality on metric measure
  spaces}}, Preprint. Available at: arXiv:1303.4382.

\bibitem{GRS_HJ}
N.~Gozlan, C.~Roberto, and P.-M. Samson, \emph{Hamilton-jacobi equations on
  metric spaces and transport entropy inequalities}, To appear in
  Rev.~Mat.~Iberoam. Available at: arXiv:1203.2783.

\bibitem{Hsu}
E.~P. Hsu, \emph{Stochastic analysis on manifolds}, Graduate studies in
  mathematics, 38, American mathematical society, Providence, RI, 2002.

\bibitem{Kend}
W.~Kendall, \emph{Nonnegative {R}icci curvature and the {B}rownian coupling
  property}, Stochastics \textbf{19} (1986), 111--129.

\bibitem{K8}
K.~Kuwada, \emph{Couplings of the {B}rownian motion via discrete apporoximation
  under lower {R}icci curvature bounds}, Probabilistic Approach to Geometry,
  Adv. Stud. Pure Math. 57, Math. Soc. Japan, 2010, pp.~273--292.

\bibitem{K9}
\bysame, \emph{Duality on gradient estimates and {W}asserstein controls}, J.
  Funct. Anal. \textbf{258} (2010), no.~11, 3758--3774.

\bibitem{K10}
\bysame, \emph{Convergence of time-inhomogeneous geodesic random walks and its
  application to coupling methods}, Ann. Probab. \textbf{40} (2012), no.~5,
  1945--1979.

\bibitem{K11}
\bysame, \emph{A probabilistic approach to the maximal diameter theorem}, Math.
  Nachr. \textbf{286} (2013), no.~4, 374--378.

\bibitem{K13}
\bysame, \emph{Gradient estimate for {M}arkov kernels, {W}asserstein control
  and {H}opf-{L}ax formula}, 
Potential Theory and Related Fields, RIMS K\^oky\^uroku Bessatsu. 
B43 (2013) 61--80. 

\bibitem{K-Phili2}
K.~Kuwada and R.~Philipowski, \emph{Coupling of {B}rownian motion and
  {P}erelman's {$\mathcal{L}$}-functional}, J. Funct. Anal. \textbf{260}
  (2011), no.~9, 2742--2766.

\bibitem{K12}
K.~Kuwada and K.-T. Sturm, \emph{Monotonicity of time-dependent transportation
  costs and coupling by reflection}, Potential Anal. \textbf{39} (2013), 
no. 3, 231--263. 

\bibitem{Led_geom-Markov}
M.~Ledoux, \emph{The geometry of {M}arkov diffusion generators}, Ann. Fac. Sci.
  Toulouse Math. (6) \textbf{9} (2000), no.~2, 305--366.

\bibitem{Lisini:2007be}
S.~Lisini, \emph{Characterization of absolutely continuous curves in
  {W}asserstein spaces}, Calc. Var. Partial Differential Equations \textbf{28}
  (2007), no.~1, 85--120.

\bibitem{Lott-Vill_AnnMath09}
J.~Lott and C.~Villani, \emph{Ricci curvature for metric-measure spaces via
  optimal transport}, Ann. Math. \textbf{169} (2009), no.~3, 903--991.

\bibitem{NP}
R. W. Neel and I. Popescu, 
\emph{A stochastic target approach to Ricci flow on surfaces},
Preprint. Available at: arXiv:1209.4150. 

\bibitem{Qian:con}
Z.-M. Qian, \emph{On conservation of probability and the {F}eller property},
  Ann. Prob. \textbf{24} (1996), no.~1, 280--292.

\bibitem{Rajala:2012jm}
T.~Rajala, \emph{{Interpolated measures with bounded density in metric spaces
  satisfying the curvature-dimension conditions of Sturm}}, J. Funct. Anal.
  \textbf{263} (2012), no.~4, 896--924.

\bibitem{Savare:2013wa}
G.~Savar{\'e}, \emph{{Self-improvement of the Bakry-{\'E}mery condition and
  Wasserstein contraction of the heat flow in RCD($K,\infty$) metric measure
  spaces}}, 
Disc. Cont. Dyn. Sist. \textbf{34} (2014), no.~4, 1641--1661. 
% Preprint. Available at: arXiv:1304.0643.

\bibitem{Sturm_Ric}
K.-T. Sturm, \emph{On the geometry of metric measure spaces. {I}}, Acta. Math.
  \textbf{196} (2006), no.~1, 65--131.

\bibitem{Sturm_Ric2}
\bysame, \emph{On the geometry of metric measure spaces. {II}}, Acta. Math.
  \textbf{196} (2006), no.~1, 133--177.

\bibitem{Topp_Lopt}
P.~Topping, \emph{{$\mathcal{L}$}-optimal transportation for {R}icci flow}, J.
  Reine Angew. Math. \textbf{636} (2009), 93--122.

\bibitem{book_Vil2}
C.~Villani, \emph{Optimal transport, old and new}, Grundlehren der
  Mathematischen Wissenschaften, vol.~338, Springer-Verlag, 2008.

\bibitem{Renes_poly}
M.-K. von Renesse, \emph{Intrinsic coupling on {R}iemannian manifolds and
  polyhedra}, Electron. J. Probab. \textbf{9} (2004), no.~14, 411--435.

\bibitem{Wang_book05}
F.-Y. Wang, \emph{Functional inequalities, {M}arkov semigroups, and spectral
  theory}, Mathematics Monograph Series 4, Science Press, Beijing, China, 2005.

\bibitem{Wang_equivCD}
\bysame, \emph{Equivalent semigroup properties for curvature-dimension
  condition}, Bull. Sci. Math. \textbf{135} (2011), no.~6--7, 803--815.

\bibitem{Wang_book14}
\bysame, \emph{Analysis for diffusion processes on Riemannian manifolds}, 
World Scientific Publishing Co. Pte. Ltd., Hackensack, NJ, 2014. 

\end{thebibliography}
\providecommand{\bysame}{\leavevmode\hbox to3em{\hrulefill}\thinspace}
\providecommand{\MR}{\relax\ifhmode\unskip\space\fi MR }
\providecommand{\MRhref}[2]{%
  \href{http://www.ams.org/mathscinet-getitem?mr=#1}{#2}
}
\providecommand{\href}[2]{#2}

\end{document}